\documentclass[11pt,times]{article}
\usepackage[applemac]{inputenc}
\usepackage{amssymb}
\usepackage{geometry}
\usepackage{amsthm}

\usepackage{lscape}

\usepackage[dvipsnames]{xcolor}

\usepackage{bm}
\usepackage{envmath}
\usepackage{subfigure}
\usepackage{graphicx}
\usepackage{natbib}
\usepackage{amsmath}

\newtheorem{theorem}{Theorem}[section]
\newtheorem{corollary}[theorem]{Corollary}
\newtheorem{lemma}[theorem]{Lemma}
\newtheorem{remark}[theorem]{Remark}

\begin{document}

\title{Multiple Hypotheses Testing For Variable Selection} 

\author{Florian Rohart}
\date{April 2012}
\maketitle
\begin{abstract}
Many methods have been developed to estimate the set of relevant variables in a sparse linear model $Y= X\beta+\epsilon$ where the dimension $p$ of $\beta$ can be much higher than the length $n$ of $Y$. Here we propose two new methods based on multiple hypotheses testing, either for ordered or non-ordered variable selection. 
Our procedures are inspired by the testing procedure proposed by \cite{Baraud:2003}. The new procedures are proved to be powerful under some conditions on the signal and their properties are non asymptotic. They gave better results in estimating the set of relevant variables than both the False Discovery Rate (FDR) and the Lasso, both in the common case $(p< n)$ and in the high-dimensional case $(p \ge n)$. 
\end{abstract}

\section{Introduction}

 Recent technologies have provided scientists with very high-dimensional data. This is especially the case in biology with high-throughput DNA/RNA chips. Unravelling the relevant variables -genes for example- underlying an observation is a well known problem in statistics and is still one of the current major challenges. Indeed, with a large number of variables there is often a desire to select a smaller subset that not only fits almost as well as the full set of variables, but also contains the most important ones for a prediction purpose. Discovering the relevant variables leads to higher prediction accuracy, an important criterion in variable selection. 
  
 Many methods have been developed to estimate the set of relevant variables in  the linear model $Y= X\beta+\epsilon$ where the dimension $p$ of $\beta$ can be much higher than the length $n$ of $Y$. Most of these methods are based on a penalized criterion.
  The mostly known is  probably the Lasso that has been presented by \cite{Tibshirani:1996}; $l^1$ penalization of the least squares estimate which shrinks to zero some irrelevant coefficients, hence an estimation of the set of relevant variables. A lot of studies have been conducted on the Lasso and many results are available; e.g. consistency \citep{Zhao:2006}, sparsity oracle inequalities \citep{Bunea:2007} and variable selection in high-dimensional graphs with the Lasso \citep{Meinshausen:2006}. The Lasso has several variants such as an adaptative Lasso \citep{Huang:2008}, a bootstrap Lasso \citep{bach:2009} or a Group Lasso \citep{Chesneau:2008}.
A $l^1$ penalization has also been used in the Sparse-PLS, which induces a limited number of variables in each PLS direction; see \citet{Tenenhaus:1998} for an introduction on PLS, and \citet{LeCao:2008} for further details on Sparse-PLS.
Other kinds of penalization have also been used, such as the Akaike Information Criterion (AIC) or the Bayesian Information Criterion (BIC), two methods based on the logarithm of the likelihood penalized by the number of variables included in the model.
Despite that the major portion of model selection methods was developed to perform in low dimension, some of them apply in the high-dimensional case. There is still some others that were actually developed to be powerful when $p$ is higher than $n$, such as the Dantzig selector \citep{Candes:2007}. Yet, a recent paper  shows that under a sparsity condition on the linear model, the Dantzig selector and the Lasso exhibit similar behavior \citep{Bickel:2009}.
Nevertheless, penalization criterion is not the only way to perform model selection.
For instance, the False Discovery Rate (FDR) procedure, developed in the context of multiple hypotheses testing by \citet{Benjamini:1995}, was used in variable selection by \citet{bunea:2006}. This procedure has been extended to high-dimensional analysis and is presently used in biology for QTL research and transcriptome analysis; a p-value is calculated for each variable $X_i$ from the regression of Y onto that variable and selection is performed through an adjusted threshold.

Most of the selection methods cited above give quite good results when $p$ is lower than $n$. However, they all have drawbacks that especially appear in a high-dimensional context. For instance, Lasso lacks stability; only small changes in the data set leads to different sets of selected variables. Moreover, the results of the Lasso, as well as its variants, depends on a penalty parameter that has to be tuned, which is surely the major drawback. \\
This paper deals with the problem of selecting the set of indices of the relevant variables in a sparse linear model when $p$ can be lower or higher than $n$. We present a new method of variable selection based on multiple hypotheses testing which is stable and free of tuning parameters.\\

We consider the regression model:
\begin{equation}
\label{0}
Y=X\beta +\epsilon,
\end{equation}
where $Y$ is the observation of length $n$, $X=(X_1,\dots,X_p)$ is the $n\times p$ matrix of $p$ variables, $\beta$ is an unknown vector of $\mathbb{R}^p$, $\epsilon$ a Gaussian vector with i.i.d. components, $\epsilon \sim \mathcal{N}_n(0,\sigma^2I_n)$ where $I_n$ is the identity matrix of $\mathbb{R}^n$, and $\sigma$ some unknown positive quantity. We assume that $X_1$ is the vector of $\mathbb{R}^n$ whose coordinates are all equal to $1/\sqrt{n}$. 
 We set $J=\{j,\beta_j\neq 0\}$ and $|J|=k_0$. We denote $\beta_J=(\beta_j)_{j\in J}$. 
Let $\mu=E(Y)=X\beta$ and $\mathbb{P}_\mu$ the distribution of $Y$ obeying to model \eqref{0}.

The aim of this paper is to estimate $J$, the set of indices of the relevant variables in \eqref{0}. We distinguish two frameworks. In a first step, we only consider ordered variable selection. We define a powerful procedure for estimating $J$ under some conditions on the signal, either when $p\le n$ or when $p> n$. These properties are non asymptotic. The procedure is a multiple hypotheses testing method based on the testing procedure developed by \citet{Baraud:2003}. 
In a second step, the variables are not assumed to be ordered. We provide a procedure to estimate $J$ when $\sigma$ is known and another procedure when $\sigma$ is unknown. The two procedures are proved to be powerful under some conditions on the signal. The properties of the procedures are also non asymptotic.  \\

Let us introduce some notations that will be used throughout this paper. Note $||s||_n^2=\sum_{i=1}^n s_i^2/n$. Set $\Pi_V$ the orthogonal projector onto $V$ for all subspace $V$.
  $\bar{F}_{D,N}(u)$ denotes the probability for a Fisher variable with $D$ and $N$ degrees of freedom to be larger than $u$.  We denote $\forall \ (x,y)\in \mathbb{R}^{n^2}$ $<x,y>_n=\sum_{i=1}^n x_iy_i/n$, $<x,y>=n<x,y>_n$ and $\forall a \in \mathbb{R}, \  \lfloor a \rfloor $ the integer part of a. \\ 
  
This paper is organized as follow, in Section \ref{ordered} we present the first procedure to estimate $J$ in the context of ordered variable selection; the non-ordered variable selection is considered in Section \ref{nonordered}. A simulation study is provided in Section \ref{result} to compare several variable selection methods. The proofs are given in Appendix \ref{proof}.

\section{Ordered variable selection}
\label{ordered}
\subsection{The low-dimensional case}
\label{ordered1}
First of all, the common case $p< n$ is considered.
The variables $(X_i)_{1\le i\le p}$ are assumed to be linearly independent and for all $i \in \{1,\dots,p\}$, $X_i$ is supposed to be normed to $1$: $\forall \ i, ||X_i||_n=1$. 

In this section we focus on ordered variables selection, which means that the set of indices of the relevant variables is supposed to be $J=\{1,\dots,k_0\}$, for some $k_0\le p$. Hence an estimation of $k_0$ gives us an estimation of $J$. This section focuses on the estimation of $k_0$.\\
Our procedure is a multiple hypotheses testing method based on the testing procedure developed by \citet{Baraud:2003}
 in the context of linear regression of $Y=f+\epsilon$ where $f$ is an unknown vector of $\mathbb{R}^n$. Let $V$ be a subspace of $\mathbb{R}^n$. They constructed a testing procedure of the null hypothesis ``$f$ belongs to $V$" against the alternative that it does not under no prior assumption on $f$. Their testing procedure is based on the choice of a collection $\{S_{m}, m\in \mathcal{M}\}$ of subspaces of $V^\bot$ and the choice of a collection of levels $\{\alpha_m,m\in \mathcal{M}\}$. They considered  for each $m\in \mathcal{M}$ a Fisher test of level $\alpha_m$ to test
\begin{equation*}
 H_0: \{f\in V\} \hspace{0.5cm} \text{against the alternative} \hspace{0.5cm} H_{1,m}: \{f\in (V+S_m)\backslash V\},
 \end{equation*}
and the null hypothesis $H_0$ is rejected if  one of the Fisher test does.\\

Our procedure consists in applying the procedure proposed by \citet{Baraud:2003} on a collection of subspaces $(V_k)_{1\le k< p}$ to test successively the null hypotheses $H_k: \{\mu \in V_k\}$, for $1\le k < p$. We stop our procedure as soon as a null hypothesis is accepted.\\
Set $\forall 1\le k < p, V_k=\text{span}(X_1,\dots,X_{k})$. With this choice of $V_k$ and as $(X_i)_{1\le i\le p}$ is a linearly independent family, we have $\text{dim}(V_k)=k,\forall 1\le k < p$.\\
Let $k$ be fixed in $\{1,\dots,p-1\}$, we define $t_{max}^k=\lfloor log_2(p-k)\rfloor$ and $\mathcal{T}_k=\{0,\dots,t^k_{max} \}$. \\
As done in \citet{Baraud:2003}, given a collection of level $\{\alpha_{k,t}, \ t\in \mathcal{T}_k\} $ and a collection of spaces $\{S_{k,t}, \ t \in \mathcal{T}_k\}$
we consider for each $t\in \mathcal{T}_k$ a Fisher test of level $\alpha_{k,t}$ to test the null hypothesis 
 \begin{equation*}
H_k:  \{\mu \in V_k\}\hspace{0.5cm} \text{ against the alternative }
\hspace{0.5cm} \{\mu \in (V_k+S_{k,t})\backslash V_k\}.
 \end{equation*}
 The null hypothesis $H_k$ is rejected if at least one of the Fisher tests does. 
The collection of levels $\{\alpha_{k,t}, \ t\in \mathcal{T}_k\} $ is calibrated in order to ensure that the final test $H_k$ is of level $\alpha$ -fixed in $]0,1[$-, and the collection $\{S_{k,t}, \ t \in \mathcal{T}_k\}$ of subspaces of $V^\bot$ is defined as follows:\\
 $\forall \  t\in \mathcal{T}_k$, 
\begin{equation}
\label{Q}
S_{k,t} =  \text{span}\left(\Pi_{V_k^\bot}(X_{k+1}),\dots, \Pi_{V_k^\bot}(X_{k+2^t})\right).
\end{equation} 
Let us introduce some notations that will be used throughout this section. For each $k \in \{1,\dots,p\}$, $t\in \mathcal{T}_k$, we set $V_{k,t}=V_k\oplus S_{k,t}$, and denote $D_{k,t}=2^t$ and $N_{k,t}=n-(k+2^t)$ the dimension of $S_{k,t}$ and $V_{k,t}^\bot$ respectively. \\

As our procedure consists in successively testing the null hypotheses $(H_k)_{1\le k < p}$ at level $\alpha$ until a null hypothesis is accepted, an estimation of $k_0$ with our procedure is $$\hat{k}=\text{inf} \{k\ge1, H_{k}\ \text{is accepted}\}.$$ The estimated set of indices of the relevant variables is then $\hat{J}=\{1,\dots,\hat{k}\}$. Note that if all the null hypotheses $(H_k)_{1\le k<p}$ are rejected, $\hat{J}=\{1,\dots,p\}$.\\
Let us recall the definition of the procedure proposed by \citet{Baraud:2003} to test the null hypothesis $H_k$.
Set: $\forall \alpha \in ]0,1[$, $\forall 1\le k<p$,
  \begin{equation}
  \label{Tk}
   T_{k,\alpha}=\underset{t\in \mathcal{T}_k}{\text{sup}} \{ \dfrac{N_{k,t} ||\Pi_{S_{k,t}} Y||_n^2}{D_{k,t} ||Y-\Pi_{V_{k,t}}Y||_n^2} -\bar{F}_{D_{k,t},N_{k,t}}^{-1}(\alpha_{k,t})    \},
\end{equation}
where $\{\alpha_{k,t}, \ t\in \mathcal{T}_k\} $ is a collection of number in ]0,1[ such that:
\begin{equation}
\forall \mu \in V_k, \hspace{1cm} \mathbb{P}_\mu (T_{k,\alpha} >0) \le \alpha.
 \end{equation}
 The null hypothesis $H_{k}$ is rejected when $T_{k,\alpha}$ is positive.\\
They chose the collection $\{\alpha_{k,t}, t\in \mathcal{T}_k\}$ in accordance with one of the two following procedures:\\
 
 P1. For all $ t\in \mathcal{T}_k, \alpha_{k,t}=\alpha_{k,n}$ where $\alpha_{k,n}$ is the $\alpha$-quantile of the random variable
 \begin{equation*}
 \underset{t\in \mathcal{T}_k}{\text{inf}} \ \bar{F}_{D_{k,t},N_{k,t}}\{\dfrac{N_{k,t} ||\Pi_{S_{k,t}} \epsilon||_n^2}{D_{k,t} ||\epsilon-\Pi_{V_{k,t}}\epsilon||_n^2}  \},
 \end{equation*}
 
 P2. The collection $\{\alpha_{k,t},t\in \mathcal{T}_k\}$ satisfies the inequality
 \begin{equation*}
\sum_{t\in \mathcal{T}_k}\alpha_{k,t}\le\alpha.
 \end{equation*}
\\
Procedure P1 gives a test $H_{k}$ of size $\alpha$ whereas procedure P2 only gives a test $H_{k}$ of level $\alpha$.
Our final multiple testing procedure, which consists in calculating successively $T_{k,\alpha}$ from $k=1$ until $T_{k,\alpha}$ is negative, is proved to be powerful. An upper bound of the probability to wrongly estimate $k_0$ is given in the following theorem.
Let us first introduce some notations. For $k=1,\dots,p-1$, for $\gamma \in ]0,1[$ and for all $t\in \mathcal{T}_k$, let $L_t=\text{log}(1/\alpha_{k,t}), L=\text{log}(2/\gamma), m_t=2\text{exp}(4L_t/N_{k,t})$, and for $u>0$ let
$$K_t(u)=1+2\sqrt{\dfrac{u}{N_{k,t}}}+2m_t\dfrac{u}{N_{k,t}},$$
$$C_1(k,t)=2.5(1+K_t(L_t)\vee m_t)\dfrac{D_{k,t}+L_t}{N_{k,t}},$$
$$C_2(k,t)=2.5\sqrt{1+K_t^2(L)}\left(1+\sqrt{\dfrac{D_{k,t}}{N_{k,t}}}\right),$$
$$C_3(k,t)=2.5\left[\left(\dfrac{m_tK_t(L)}{2}\right)\vee 5\right]\left(1+2\dfrac{D_{k,t}}{N_{k,t}}\right),$$

\begin{theorem}
\label{th_1}
Let $Y$ obey to Model \eqref{0}. We assume that $p< n$ and that  $(X_i)_{1\le i\le p}$ is a linearly independent family. We denote by $J$ the set $\{j, \beta_j\neq 0\}=\{1,\dots,k_0\}$. 
Let $\gamma \ \text{and } \alpha$ be fixed in $]0,1[$ .\\
The testing procedure estimates $k_0$ by $\hat{k}=\text{inf} \{k\ge1, T_{k,\alpha}\le 0\}$, where $T_{k,\alpha}$ is defined by \eqref{Tk}. Let $\{\alpha_{k,t}, t\in \mathcal{T}_k\}$ be defined according to the procedure P1 or P2. \\
The following inequality holds for all $\mu\in \mathbb{R}^n$ and for all $k_0\le p$:
\begin{equation}
\label{sup}
\mathbb{P}_{\mu}(\hat{k}> k_0) \le \alpha.
\end{equation}

If \ $\forall k \le k_0-1$ the condition $(R_k)$ holds\\
$(R_k) : \exists \ t \in \mathcal{T}_k/ $
\begin{eqnarray*}
\left|\left|\Pi_{ S_{k,t}}(\mu)\right|\right|_n^2 &\ge& C_1(k,t)\left|\left|\Pi_{V_{k,t}^{\bot}}(\mu)\right|\right|_n^2\\
&& +\dfrac{\sigma^2}{n} \left[ C_2(k,t) \sqrt{2^t \ log\left(\dfrac{2k_0}{\alpha_{k,t} \gamma}\right)} + C_3(k,t) log\left(\dfrac{2k_0}{\alpha_{k,t} \gamma}\right)\right]
\end{eqnarray*}
then 
\begin{equation}
\label{inf}
\mathbb{P}_{\mu}(\hat{k}< k_0) \le  \gamma,
\end{equation}
which implies that
\begin{equation}
\label{non}
\mathbb{P}_{\mu}(\hat{k}\neq k_0) \le  \gamma + \alpha.
\end{equation}
\end{theorem}

This result is derived from the result on the power of the multiple testing procedure proposed by \citet{Baraud:2003}. It is important to note that Theorem \ref{th_1} is non asymptotic. \\

{\bf{Comments}}
\begin{enumerate}
\item
\label{rem1}
As mentioned in \citet{Baraud:2003}, for $k$ fixed, $C_1(k,t)$, $C_2(k,t)$ and $C_3(k,t)$ behave like constants if the following conditions are verified:

For all $t\in \mathcal{T}_k, \alpha_{k,t}\ge \text{exp}(-N_{k,t}/10), \  \gamma\ge 2\text{exp}(-N_{k,t}/21)$ and the ratio $\dfrac{D_{k,t}+L_{k,t}}{ N_{k,t}}$ remains bounded.\\\\
Under these conditions, the following inequalities hold:
$$C_1(k,t)\le 10 \dfrac{D_{k,t}+\text{log}(1/\alpha_{k,t})}{N_{k,t}},$$ $$C_2(k,t)\le5\left(1+\sqrt{\dfrac{D_{k,t}}{N_{k,t}}}\right),$$
$$C_3(k,t)\le12.5\left(1+2\dfrac{D_{k,t}}{N_{k,t}}\right).$$

\item
We say that $\mu$ satisfies condition $(R)$ if $\forall k \le k_0-1$, $(R_k)$ holds. According to Theorem \ref{th_1}, our procedure is powerful under the condition $(R)$. A condition on the coefficients $\beta_J$ underlies in $(R)$ since the projection of $Y$ onto a space spanned by a subset of the family $(X_i)_{1\le i\le p}$ depends both on $\beta$ and on the matrix $X$. These conditions on $\beta_J$ explicitly appear when $(X_i)_{1\le i\le p}$ is an orthonormal family. Assuming that $(X_i)_{1\le i\le p}$ is an orthonormal family, \eqref{0} becomes:
\begin{equation}
\label{1}
Y=\underbrace{X_1\beta_1+..+X_k\beta_k}_{\in V_k}+\underbrace{X_{k+1}\beta_{k+1}+...+X_p\beta_p}_{\in V_k^{\bot}}+\epsilon.
\end{equation}
With the new decomposition \eqref{1}, the projection of $Y$ on any subspace $S_{k,t}$ only depends on the coefficients $(\beta_j)_{j\ge k+1}$.
Thus the condition $(R_k)$ can be written in a more explicit form, involving the coefficients $(\beta_j)_{1\le j \le p}$. Namely, $(R_k)$ is equivalent to:
$ \exists \ t\in\mathcal{T}_k /$
\begin{eqnarray*} 
\beta_{k+1}^2+..+\beta_{k+2^t}^2 &\ge& \ C_1(k,t)\sum_{j=k+2^t+1}^p \beta^2_j\\
&& +\dfrac{\sigma^2}{n} \left[ C_2(k,t) \sqrt{2^tlog\left(\dfrac{2k_0}{\alpha_{k,t} \gamma}\right)}+C_3(k,t) log\left(\dfrac{2k_0}{\alpha_{k,t} \gamma}\right)\right] .
\end{eqnarray*}
When $k<k_0$, the coefficients $\beta_{k+1},\dots,\beta_{k_0}$ are not equal to $0$. If for some $t\in \mathcal{T}_k$, the sum $\beta_{k+1}^2+\dots+\beta_{k+2^t}^2$ is large enough (namely larger than the right hand of the above equation), then the test will be powerful and the hypotheses $H_k$ will be rejected with high probability.\\
\end{enumerate}

\subsection{The high-dimensional case}
\label{highdimsec}
After pointing out the properties of the procedure in the case $p< n$, let us discuss the high-dimensional case, i.e. $p\ge n$. The family $(X_i)_{1\le i\le p}$ can no longer be a linearly independent family. We have to make the assumption that the decomposition of $\mu$ is unique, i.e $\exists ! J\subset\{1,\dots,p\} /\mu=\sum_{j\in J}X_j\beta_j$. \\
 Let recall that in this section we only considered ordered variable selection. We then assume $J=\{1,\dots,k_0\}$, hence $|J|=k_0$.
 
Let denote by $a_s$ the dimension of the linear span of $(X_1,\dots,X_s)$ for $1\le s \le p$, note that $a_p\le n$. 
We define, for all $ k< a_p-1$, $t^k_{max}=\lfloor log_2(a_p-k-1)\rfloor$ and $\mathcal{T}_k=\{0,\dots,t^k_{max} \}$. Denote $V_{k}=\text{span}(X_1,\dots,X_{s_k})$ where $s_k$ is defined by $s_k=\text{inf}\{s: a_s=k\}$. 
As in the previous section, we test successively the hypotheses $H_k: \{\mu \in V_k\}$ for all $1\le k<a_p-1$ until a null hypothesis is accepted.\\
 In order to test $H_k$, we consider the collection of alternatives $\{\mu \in V_k+S_{k,t}\}$ for $t\in \mathcal{T}_k$ where $S_{k,t}=\text{span}\left(\Pi_{V_k^\bot}(X_{s_k+1},\dots,\Pi_{V_k^\bot}(X_{s_k+q_{k,t}})\right)$ with $q_{k,t}$ defined by 
 $q_{k,t}=\text{inf}\{q:\text{dim}\left(Q_{k,q}\right)=2^t\}$, 
where for all $q>0$,
$Q_{k,q}=\text{span}\left(\Pi_{V_k^\bot}(X_{s_k+1}),\dots,\Pi_{V_k^\bot}(X_{s_k+q})\right)$.

The condition $(R_k)$ in Theorem \ref{th_1} gives no restriction on the growth of $p$. Thus Theorem \ref{th_1} applies with the new notations for any $p\ge n$, but for $k_0<n$. This is no strong restriction since $k_0>n$ means that no sufficient observations are available  to estimate $k_0$, whatever the method.

Results from a simulation study in Section \ref{result} will show the power of our procedure; either when $p< n$ or when $p\ge n$.

\section{Non-ordered variable selection}
\label{nonordered}

In Section \ref{ordered} we defined a procedure based on multiple hypotheses testing in order to estimate $J$, the set of indices of the relevant variables of a sparse linear model \eqref{0}. As we considered ordered variable selection, the estimation of $J=\{1,\dots,k_0\} $ was reduced to the estimation of $k_0$. The present section is dedicated to non-ordered variable selection, so  $J$ is not necessarily equal to $\{1,\dots,k_0\}$. 
We define here a general two-step procedure to estimate $J$; the first step orders the variables and the second performs multiple testing. After the first step of the general procedure, the ordered variables will be denoted as $X_{(1)},\dots,X_{(p)}$, where $X_{(1)}=X_1$. \\

The first step of our procedure consists in ordering the variables. In this paper, we proposed two ways to order  $(X_i)_{2\le i\le p}$ taking into account the observation $Y$.\\
\begin{enumerate}
\item Variables ordered by increasing p-values: when $p< n$, a p-value is calculated for each variable from the test of nullity of the coefficient associated to this variable and the variables are sorted by increasing p-values. When $p\ge n$, a p-value is calculated for each variable using the decomposition of $Y$ onto that variable.\\ 

\item The second method that we propose orders the variables with the Bolasso technique, introduced by \citet{bach:2009}. It is a bootstrapped version of the Lasso which improves its stability: several independent bootstrap samples are generated and the Lasso is performed on each of them. This approach is proved to make the irrelevant variables asymptotically disappear. 
 A modification is applied to the Bolasso to adapt it to non asymptotic analysis.
An appearance frequency is calculated for each variable $X_i$ by counting the number of times the 
variable $X_i$ is selected over the bootstrap samples. A high frequency denotes a good prediction 
ability of the variable $X_i$, at a given penalty. To avoid the use of a penalty, we set the first ordered variable of the family $(X_i)_{2\le i\le p}$ to be the first one to reach the frequency $1$ from a decreasing penalty; and so on for the other variables. We proceed by dichotomy to order the variables. 
\
\end{enumerate}

The first method has been considered since it is often used in practice, in particular in the False Discovery Rate procedure \citep{Benjamini:1995}.
It is the one requiring less computational time, but as shown in Section \ref{result}, the Bolasso technique gives better results. Indeed, as we will see throughout this section, the first step is crucial; the ability to estimate $J$ with this procedure depends on the ability to get the relevant variables in the first places. \\

From now on, we assume that the decomposition of $\mu$ is unique, i.e $\exists ! J\subset\{1,\dots,p\} /\mu=\sum_{j\in J}X_j\beta_j$. We have $|J|=k_0$. A variable $X_j$ is said to be relevant when $\beta_j\neq 0$.
We introduce here an event that will be useful in the following of this section:
\begin{equation}
\label{Ak}
A_k=\{ \{(1),\dots,(k)\} = J  \}.
\end{equation}

On the event $A_{k}$,  the set of the $k$ first ordered variables corresponds to the set $J$ of the relevant variables.
The second step of the general procedure consists in testing successively the null hypothesis:
\begin{equation}
\label{Hk}
\hat{H}_k:\hspace{0.2cm} \{\mu \in \text{span}(X_{(1)},\dots,X_{(k)})\} \text{ against the alternative }\\\text{ that it does not. }
\end{equation}
The procedure stops when the null hypothesis is accepted: $$\mathring{k}=\text{inf}\{k\ge1, \hat{H}_k \text{ is accepted}\}.$$
We estimate the set $J$ of relevant variables by $$\hat{J}=\{(1),\dots,(\mathring{k})\}.$$

Note that this is not a simple generalization of the procedure proposed in Section \ref{ordered} since $\text{span}(X_{(1)},\dots,X_{(k)})$ are random spaces depending on the observation $Y$ which have been used in the first step to order the family $(X_i)_{2\le i\le p}$. The same observation $Y$ will be used in the second step to perform the multiple testing procedure.

For the sake of understanding, we first deal with the case where $\sigma$ is known in order to propose a multiple testing procedure. 

\subsection{The low-dimensional case with known variance}
\label{sigyes}
In this section, we define a procedure called Procedure `A' under the assumption that the variance $\sigma^2$ is known. 
Assume that $(X_i)_{1\le i \le p}$ is a linearly independent family and that the first step of Procedure `A' has already been done: variables have been ordered. The second step is a testing procedure that will be described in the following. 
As in the previous section, we test successively the null hypotheses $\hat{H}_k$ for $1\le k < p$ until a null hypothesis is accepted.

Let us adapt the notation of Section \ref{ordered1} to this section: we first recall that $\forall \ 1\le k\le p-1$, $t^k_{max}=\lfloor log_2(p-k)\rfloor$, $\mathcal{T}_k=\{0,\dots,t^k_{max}\}$. We define  $V_{(k)}=\text{span}(X_{(1)},\dots,X_{(k)})$ and  $\forall \  t\in\mathcal{T}_k,$\\
$S_{(k),(t)} =  \text{span}\left(\Pi_{V_{(k)}^\bot}(X_{(k+1)}),\dots, \Pi_{V_{(k)}^\bot}(X_{(k+2^t)})\right)$.\\
With the definition of $S_{(k),(t)}$, we have $\text{dim}(S_{(k),(t)})=D_{k,t}=2^t$.
Let us denote $ V_{(k),(t)}= V_{(k)} \oplus S_{(k),(t)}$. \\
For all $ t\in\mathcal{T}_k$, our aim is to test
\begin{equation}
\label{h1}
\hat H_k:  \{\mu \in V_{(k)}\}\hspace{0.5cm} \text{ against the alternative }\\\hspace{0.5cm} \{\mu \in (V_{(k)}+S_{(k),(t)})\backslash V_{(k)}\}.
\end{equation}

Since the variance is assumed to be known, we introduce for all $1\le k <p$ and for all $ t \in \mathcal{T}_k$,\\
$${U}_{k,t}=\dfrac{||\Pi_{S_{(k),(t)}} Y||_n^2}{\sigma^2}.$$
We introduce a multiple testing procedure that relies on the statistics $\{U_{k,t}, t\in \mathcal{T}_k\}$.

Since the spaces $\{S_{(k),(t)}, t\in \mathcal{T}_k\}$ are random and depend on $Y$ as mentioned before, we first provide a stochastic upper bound for the statistics $U_{k,t}$ in order to define the multiple testing procedure.\\

Let $\epsilon'\sim \mathcal{N}_n(0,\sigma^2I_n)$. For all $1\le k \le p$, we define a permutation $\sigma_1^k$ of $\{1,\dots,p\}$:\\
$\sigma_1^k(j)=(j)$ for all $j\in\{1,\dots,k\}$.\\
For $j\in\{k+1,\dots,p\}$, set $X_i^{(j)}=\Pi_{\text{span}(X_{\sigma_1^k(1)},\dots, X_{\sigma_1^k(j-1)})^\bot}(X_i)$ for all $1\le i\le p$ and define
$\sigma_1^{k}(j)=\underset{i\in\{1,\dots,p\}}{\text{argmax}}\left|\left|\Pi_{\text{span}(X_i^{(j)})} (\epsilon')\right|\right|_n^2$.\\

Set for all $1\le k <p$ and for all $ t \in \mathcal{T}_k$,\\
$$ U^1_{k,t}=\dfrac{||\Pi_{S_{(k),\sigma_1^k(t)}} \epsilon'||_n^2}{\sigma^2},$$
where $S_{(k),\sigma_1^k(t)}=\text{span}\left(\Pi_{V_{(k)}^\bot}(X_{\sigma_1^k(k+1)}),\dots,\Pi_{V_{(k)}^\bot}( X_{\sigma_1^k(k+2^t)})\right).$

Note that the distribution of $U^1_{k,t}$ only depends on the design matrix $X$, and can therefore be simulated for a given matrix $X$.

\begin{lemma}
\label{over12}
Let $1\le k < p$ and $t\in \mathcal{T}_k$.\\
We define $A_k=\{\{X_{(1)},\dots,X_{(k)}\}=\{X_j,j\in J\}\}$. For all $x>0$ we have
$$\mathbb{P}\left((U_{k,t}>x)\cap A_k\right)\le \mathbb{P}\left(U_{k,t}^1>x\right).$$
\end{lemma}

Let $\overline{U^1_{k,t}}(u)$ denote  the probability for the statistic $U^1_{k,t}$ to be larger than $u$.\\
Set $\forall \alpha \in ]0,1[$, $\forall 1\le k<p$,
\begin{equation}
\label{M1}
M_{k,\alpha}=\underset{t\in \mathcal{T}_k}{\text{sup}} \{U_{k,t} -\overline{U^1_{k,t}}^{-1}(\alpha_{k,t})    \},
\end{equation}
where $\{\alpha_{k,t}, \ t\in \mathcal{T}_k\} $ is a collection of number in ]0,1[ chosen in accordance to the following procedure:\\\\
 P3. For all $ t\in \mathcal{T}_k, \alpha_{k,t}=\alpha_{k,n}$ where $\alpha_{k,n}$ is the $\alpha$-quantile of the random variable
 \begin{equation*}
 \underset{t\in \mathcal{T}_k}{\text{inf}} \ \overline{U^1_{k,t}}\{U^1_{k,t}\}.
 \end{equation*}
The null hypothesis $\hat{H}_k$ is rejected when $M_{k,\alpha}$ is positive. In fact, the second step of the procedure `A' is to calculate $M_{k,\alpha}$ from $k=1$ until $M_{k,\alpha}$ is negative. The calculation of the collection $\{\alpha_{k,t}, \ t\in \mathcal{T}_k\} $ with the procedure P3 ensures that $\mathbb{P}\left((M_{k,\alpha}>0)\cap A_k\right)\le \alpha$.

In summary, the two-step procedure `A' when $\sigma$ is known consists in ordering the $p$ variables and then estimating $J$ by $\hat{J}=\{(1),\dots,(\mathring{k}_A)\}$ where $\mathring{k}_A=\text{inf}\{k\ge1, M_{k,\alpha}\le0\}$. Note that if all the null hypotheses $(\hat H_k)_{1\le k<p}$ are rejected, then $\hat{J}=\{1,\dots,p\}$.\\
The testing procedure `A' is proved to be powerful and we give an upper bound of the probability to wrongly estimate $J$ in the next theorem.

\begin{theorem}
\label{th_2}
Let $Y$ obey to Model \eqref{0}. We assume that $p< n$ and that the family $(X_i)_{1\le i\le p}$ is linearly independent. We denote by $J$ the set $\{j, \beta_j\neq 0\}$ and by $k_0$ its cardinality. 
Let $\alpha$ and $\gamma$ be fixed in $]0,1[$.\\
The procedure `A' estimates $J$ by $\hat{J}=\{(1),\dots,(\mathring{k}_A)\}$ where $\mathring{k}_A=\text{inf}\{k\ge1, M_{k,\alpha}\le0\}$, where $M_{k,\alpha}$ is defined by \eqref{M1} and $\{\alpha_{k,t}, t\in \mathcal{T}_k\}$ is defined according to the procedure P3. \\\\
We consider the condition $(R_{2,k})$ for $k< k_0$ stated as\\
$(R_{2,k}):\ \exists t\le log_2(k_0-k)$ such that
\begin{eqnarray*}
\dfrac{1}{2\sigma^2} \text{inf}\{|| \Pi_{S}\mu||_n^2,S\in B_{2^t}\}&\ge&\dfrac{2^t}{n} \left[10 +4log\left(\dfrac{(p-k)k_0}{2^{2t}}\right)\right] \\
&&+\dfrac{2}{n}\left[\sqrt{2^{t+1} log\left(\dfrac{k_0|\mathcal{T}_k|}{\gamma\alpha}\right) }+log\left(\dfrac{k_0|\mathcal{T}_k|}{\gamma\alpha}\right) \right],
\end{eqnarray*}
where $\forall d\le k_0,  B_{d}=\{\text{span}(X_I),I\subset J, |I|=d\}$ and $|\mathcal{T}_k|=\lfloor log_2(p-k)\rfloor+1$.\\\\
If  \ $\forall k \le k_0-1$ the condition $(R_{2,k})$ holds, then 
\begin{equation}
\label{non1}
\mathbb{P}_{\mu}(\hat{J}\neq J) \le   \gamma+\alpha + \delta,
\end{equation}
where $\delta=\mathbb{P}_{\mu}(A_{k_0}^c)= P_{\mu}(\exists \ j\le k_0 / \beta_{(j)}=0) $. 
\end{theorem}

This theorem is non asymptotic and its result differs from Theorem \ref{th_1} on the right part of \eqref{non1}. Indeed, the weight of the first step of the procedure, which lies in $\delta$, was not involved in Section \ref{ordered} since we considered ordered variable selection. Theorem \ref{th_2} shows that the first step is essential in the two-step procedure `A', which is easily understandable since there is no chance of having $J=\hat{J}$ if the event $A_{k_0}$ does not occur at the end of the first step of procedure `A'. Moreover, the condition $(R_{2,k})$ is also more restrictive than the condition $(R_k)$ which appeared in Theorem \ref{th_1}.\\

When the family $(X_i)_{2\le i\le p}$ is orthonormal, $U_{k,t}$ can be stochastically upper bounded by a statistic that does not depend on $(X_i)_{1\le i\le p}$ nor on the first step of the procedure. \\
Let $D>0$ and $W_1,\dots,W_{D}$ be $D$ i.i.d. standard Gaussian variables ordered as $|W_{(1)}|\ge\dots\ge|W_{(D)}|$.\\ We define: 
$\forall d=1,\dots,D$, 
\begin{equation}
\label{Zd}
Z_{d,D}=\sum_{j=1}^d W_{(j)}^2.
\end{equation}
Let $\bar{Z}_{d,D}(u)$ denote  the probability for the statistic $Z_{d,D}$ to be larger than $u$.\\
A multiple testing procedure can be derived from procedure `A' with this upper bound. 

\begin{lemma}
\label{over1}
Let $1\le k < p$ and $t\in \mathcal{T}_k$.\\
We define $A_k=\{\{X_{(1)},\dots,X_{(k)}\}=\{X_j,j\in J\}\}$. For all $x>0$, we have

$$\mathbb{P}\left((U_{k,t}>x)\cap A_k\right)\le \mathbb{P}\left(Z_{D_{k,t},p-k}/n>x\right).$$
\end{lemma}
Set $\forall \alpha \in ]0,1[$, $\forall 1\le k<p$,
\begin{equation}
\label{Mkt}
{M}^1_{k,\alpha}=\underset{t\in \mathcal{T}_k}{\text{sup}} \{U_{k,t} -\bar{Z}_{D_{k,t},p-k}^{-1}(\alpha_{k,t})/n    \},
\end{equation}
where $\{\alpha_{k,t}, \ t\in \mathcal{T}_k\} $ is a collection of number in ]0,1[ chosen in accordance to the following procedure:\\\\
 P4. For all $ t\in \mathcal{T}_k, \alpha_{k,t}=\alpha_{k,n}$ where $\alpha_{k,n}$ is the $\alpha$-quantile of the random variable
 \begin{equation*}
 \underset{t\in \mathcal{T}_k}{\text{inf}} \ \bar{Z}_{D_{k,t},p-k}\{Z_{D_{k,t},p-k}\}.
 \end{equation*}
The null hypothesis $\hat{H}_k$ is rejected when $M^1_{k,\alpha}$ is positive. 
The major benefit of Procedure `A' when the family $(X_i)_{2\le i \le p}$ is orthonormal is that the upper bound of the statistics $U_{k,t}$ in Lemma \ref{over1} does not depend on the family $(X_i)_{1 \le i \le p}$ nor on the order on that family. Thus the collection $\{\alpha_{k,t}, t\in \mathcal{T}_k\}$ defined by the procedure P4 only depends on $k$ and $t$, with $p$ and $n$ fixed.\\

In the particular case where $(X_i)_{2\le i\le p}$ is an orthonormal family, we obtain the following corollary of Theorem \ref{th_2}, which is more explicit. 
 \begin{corollary}
\label{cor_2}
Let $Y$ obey to model \eqref{0}. We assume that $p< n$ and that $(X_i)_{2\le i\le p}$ is an orthonormal family. We denote by $J$ the set $\{j, \beta_j\neq 0\}$ and by $k_0$ its cardinality. 
Let $\alpha$ and $\gamma$ be fixed in $]0,1[$.\\
The procedure estimates $J$ by $\hat{J}=\{(1),\dots,(\mathring{k}_{Abis})\}$ where $\mathring{k}_{Abis}=\text{inf}\{k\ge1, M^1_{k,\alpha}\le0\}$, where $M^1_{k,\alpha}$ is defined by \eqref{Mkt} and $\{\alpha_{k,t}, t\in \mathcal{T}_k\}$ is defined according to the procedure P4. \\\\
We consider the condition $(R_{2bis,k})$ for $k<k_0$ stated as \\
$(R_{2bis,k}): \exists \ t \le log_2(k_0-k) $ such that
\begin{eqnarray*}
\dfrac{1}{2\sigma^2}\sum_{j=1}^{2^t} \beta_{\sigma_2(j)}^2 &\ge&   \dfrac{2^t}{n} \left[10 +4log\left(\dfrac{(p-k)k_0}{2^{2t}}\right)\right]\\
&&+\dfrac{2}{n}\left[\sqrt{2^{t+1} log\left(\dfrac{k_0|\mathcal{T}_k|}{\gamma\alpha}\right) }+ log\left(\dfrac{k_0|\mathcal{T}_k|}{\gamma\alpha}\right) \right],
\end{eqnarray*}
where $\sigma_2$ is defined by $|\beta_{\sigma_2(1)}|\le \dots\le |\beta_{\sigma_2(k_0)}|$ and $|\mathcal{T}_k|=\lfloor log_2(p-k)\rfloor+1$.\\\\
If  \ $\forall k \le k_0-1$ the condition $(R_{2bis,k})$ holds, then 
\begin{equation}
\label{non1bis}
\mathbb{P}_{\mu}(\hat{J}\neq J) \le   \gamma+\alpha + \delta,
\end{equation}
where $\delta=\mathbb{P}_{\mu}(A_{k_0}^c)= P_{\mu}(\exists \ j\le k_0 / \beta_{(j)}=0) $. 
\end{corollary}

\subsection{The low-dimensional case with unknown variance}
\label{signot}

In this section, we define a procedure `B' under the assumption that the variance $\sigma^2$ is unknown. Assume that $(X_i)_{1\le i \le p}$ is a linearly independent family and that the first step of the procedure `B' has already been done: variables have been ordered.\\
In this section, the notations of Section \ref{sigyes} are used:
$\forall \ 1\le k\le p-1$, $t^k_{max}=\lfloor log_2(p-k)\rfloor$, $\mathcal{T}_k=\{0,\dots,t^k_{max}\}$.\\ We define  $V_{(k)}=\text{span}(X_{(1)},\dots,X_{(k)})$ and  $\forall \  t\in\mathcal{T}_k,$
$S_{(k),(t)} =  \text{span}\left(\Pi_{V_{(k)}^\bot}(X_{(k+1)}),\dots, \Pi_{V_{(k)}^\bot}(X_{(k+2^t)})\right)$.\\
Denote for each $k \in \{1,\dots,p\}$, $t\in \mathcal{T}_k$, $V_{(k),(t)}=V_{(k)}\oplus S_{(k),(t)}$, and denote $D_{k,t}=2^t$ and $N_{k,t}=n-(k+2^t)$ the dimension of $S_{(k),(t)}$ and $V_{(k),(t)}^\bot$ respectively. 

Since the variance is assumed to be unknown, we introduce for all $1\le k <p$ and for all $ t \in \mathcal{T}_k$,
 $$\tilde{U}_{D_{k,t},N_{k,t}}=\dfrac{N_{k,t} ||\Pi_{S_{(k),(t)}} Y||_n^2}{D_{k,t} ||Y-\Pi_{V_{(k),(t)}}Y||_n^2}.$$

In order to test the null hypothesis $\hat H_k$ defined by \eqref{Hk}, we introduce a multiple testing procedure  which relies this time on the statistics $\{\tilde{U}_{D_{k,t},N_{k,t}}, t\in \mathcal{T}_k\}$.

As in Section \ref{sigyes}, we first provide a stochastic upper bound for the statistics $\tilde{U}_{D_{k,t},N_{k,t}}$ in order to define the multiple testing procedure. 

Let $\epsilon'\sim \mathcal{N}_n(0,\sigma^2I_n)$. 
Set for all $1\le k <p$ and for all $ t \in \mathcal{T}_k$,
$$\Upsilon_{k,t}=\dfrac{N_{k,t} ||\Pi_{S_{(k),\sigma_1^k(t)}} \epsilon'||_n^2}{D_{k,t} ||\epsilon'-\Pi_{V_{(k),\sigma_1^k(t)}}\epsilon'||_n^2},$$\\
where $S_{(k),\sigma_1^k(t)}=\text{span}\left(\Pi_{V_{(k)}^\bot}(X_{\sigma_1^k(k+1)}),\dots, \Pi_{V_{(k)}^\bot}(X_{\sigma_1^k(k+2^t)})\right)$, $V_{(k),\sigma_1^k(t)}=S_{(k),\sigma_1^k(t)}\oplus V_{(k)}$ and the permutation $\sigma_1^k$ is defined as in Section \ref{sigyes}.

\begin{lemma}
\label{over22}
Let $1\le k < p$ and $t\in \mathcal{T}_k$.\\
We define $A_k=\{\{X_{(1)},\dots,X_{(k)}\}=\{X_j,j\in J\}\}$. For all $x>0$ we have

$$\mathbb{P}\left((\tilde{U}_{D_{k,t},N_{k,t}}>x)\cap A_k\right)\le \mathbb{P}\left(\Upsilon_{k,t}>x\right).$$
\end{lemma}

Let $\bar{\Upsilon}_{k,t}(u)$ denote the probability for the statistic $\Upsilon_{k,t}$ to be larger than $u$.\\
Set $\forall \alpha \in ]0,1[$, $\forall 1\le k <p$,
\begin{equation}
\label{M}
\hat{M}_{k,\alpha}=\underset{t\in \mathcal{T}_k}{\text{sup}} \{\tilde{U}_{D_{k,t},N_{k,t}} -\bar{\Upsilon}_{k,t}^{-1}(\alpha_{k,t})    \},
\end{equation}
where $\{\alpha_{k,t}, \ t\in \mathcal{T}_k\} $ is a collection of number in ]0,1[ chosen    in accordance to the following procedure:\\\\
 P5. For all $ t\in \mathcal{T}_k, \alpha_{k,t}=\alpha_{k,n}$ where $\alpha_{k,n}$ is the $\alpha$-quantile of the random variable
 \begin{equation*}
 \underset{t\in \mathcal{T}_k}{\text{inf}} \ \bar{\Upsilon}_{k,t} \{\Upsilon_{k,t} \},
 \end{equation*}
The null hypothesis $\hat{H}_k$ is rejected when $\hat{M}_{k,\alpha}$ is positive. In fact, the second step of the procedure `B' is to calculate $\hat{M}_{k,\alpha}$ from $k=1$ until $\hat{M}_{k,\alpha}$ is negative. 
The calculation of the collection $\{\alpha_{k,t}, \ t\in \mathcal{T}_k\} $ with the procedure P5 ensures that $\mathbb{P}\left((\hat M_{k,\alpha}>0)\cap A_k\right)\le \alpha$.

In summary, this two-step procedure `B' when $\sigma$ is unknown consists in ordering the $p$ variables and then estimating $J$ by $\hat{J}=\{(1),\dots,(\mathring{k}_B)\}$ where $\mathring{k}_B=\text{inf}\{k\ge1, \hat{M}_{k,\alpha}\le0\}$.
The procedure is proved to be powerful in the next theorem; we give an upper bound of the probability to wrongly estimate $J$ under some conditions on the signal.
Let us introduce some notations that will be used in the following theorem. We set\\ $L_t=\text{log}(|\mathcal{T}_k|/\alpha)$,
$m_t=\text{exp}(4L_t/N_{k,t})$, $ m_p=\text{exp}\left(\dfrac{4D_{k,t}}{N_{k,t}}log\left(\dfrac{e(p-k)}{D_{k,t}}\right)\right)$,
$M=2m_tm_p$. Denote $\Lambda_1(k,t)=\sqrt{1+\dfrac{D_{k,t}}{N_{k,t}}}$, $\Lambda_2(k,t)=\left(1+2\dfrac{D_{k,t}}{N_{k,t}}\right)M$ and $\Lambda_3(k,t)=2\Lambda_1(k,t)+\Lambda_2(k,t)$.

\begin{theorem}
\label{th_3}
Let $Y$ obey to model \eqref{0}. We assume that $p< n$ and that the family $(X_i)_{1\le i\le p}$ is linearly independent. We define by $J$ the set $\{j, \beta_j\neq 0\}$ and by $k_0$ its cardinality. 
Let $\alpha$ and $\gamma$ be fixed in $]0,1[$.\\
The procedure `B' estimates $J$ by $\hat{J}=\{(1),\dots,(\mathring{k}_B)\}$ where $\mathring{k}_B=\text{inf}\{k\ge1, \hat{M}_{k,\alpha}\le0\}$, where $\hat{M}_{k,\alpha}$ is defined by \eqref{M} and $\{\alpha_{k,t}, t\in \mathcal{T}_k\}$ is defined according to the procedure P5. \\\\
We consider the condition $(R_{3,k})$  for $k< k_0$ stated as \\$(R_{3,k}):\ \exists t\le log_2(k_0-k)$ such that\\
\begin{eqnarray*}
\dfrac{1}{2} \text{inf}\{|| \Pi_{S}\mu||_n^2,S\in B_{2^t}\}\ &
\ge &\label{right}\dfrac{A(k,t)}{N_{k,t}}\left[||\mu||^2_n+\sigma^2\left(2+\dfrac{3}{n}log\left(\dfrac{2k_0}{\gamma}\right)\right)\right]\\
&  +&\dfrac{\sigma^2}{n}\left[2^t\left(6 +4log\left(\dfrac{k_0}{2^t}\right)\right)+3log\left(\dfrac{2k_0}{\gamma}\right)\right],
  \end{eqnarray*}
 where 
\begin{eqnarray*} 
A(k,t) &=& 2^t\left[ 2+\dfrac{2^{t}}{N_{k,t}}+\Lambda_3(k,t)log\left(\dfrac{e(p-k)}{2^t}\right)\right]\\
& +&\left(1+\Lambda_2(k,t)\right)log\left(\dfrac{log_2(p-k)+1}{\alpha}\right),
\end{eqnarray*}
and $\forall d\le k_0,  B_{d}=\{\text{span}(X_I),I\subset J, |I|=d\}$.\\\\
If  \ $\forall k \le k_0-1$ the condition $(R_{3,k})$ holds, then 
\begin{equation}
\label{non3}
\mathbb{P}_{\mu}(\hat{J}\neq J) \le   \gamma+\alpha + \delta,
\end{equation}
where $\delta=\mathbb{P}_{\mu}(A_{k_0}^c)= P_{\mu}(\exists \ j\le k_0 / \beta_{(j)}=0) $. 
\end{theorem}

This theorem is non asymptotic and shows that the testing procedure `B' is powerful under some conditions on the signal. As for Theorem \ref{th_2} of Section \ref{sigyes}, the first step of the procedure -the ordering of the variables- has an important part in Theorem \ref{th_3}. 
A simulation study in Section \ref{result} will show that this testing procedure combined with a good way to order variables -in order to minimize $\delta$- performs well.\\
\begin{remark}
\label{rem3}
The condition $(R_{3,k})$ can be simplified under the assumption that $2^t\le (n-k)/2$ and $log(p-k)>1$. Indeed, in this case, the right hand in condition $(R_{3,k})$ is upper bounded by 
\begin{equation}
\label{simp}
C(||\mu||_n,\gamma,\alpha,\sigma)2^t\left[\dfrac{log(p-k)}{N_{k,t}}+\dfrac{log(k_0)}{n}\right],
\end{equation} where $C(||\mu||_n,\gamma,\alpha,\sigma)$ is a constant depending on $||\mu||_n,\gamma$, $\alpha$ and $\sigma$.
\end{remark}

\begin{remark}
\label{ortho}
When $(X_i)_{2\le i\le p}$ is an orthonormal family, the condition $(R_{3,k})$ of Theorem \ref{th_3} can be rewritten in a more explicit way. The new condition $(R_{3bis,k})$ obtained in this case is the following:\\
$(R_{3bis,k}):\ \exists t\le log_2(k_0-k)$ such that\\
\begin{eqnarray*}
\dfrac{1}{2\sigma^2}\sum_{j=1}^{2^t} \beta_{\sigma_2(j)}^2 &\ge&\dfrac{A(k,t)}{N_{k,t}}
\left[\sum_{j=k+2^t}^{j=k_0}\beta^2_{\sigma_2(j)}+\sigma^2\left(2+\dfrac{3}{n}log\left(\dfrac{2k_0}{\gamma}\right)\right)\right]\\
&  +&\dfrac{\sigma^2}{n}\left[2^t\left(6 +4log\left(\dfrac{k_0}{2^t}\right)\right)+3log\left(\dfrac{2k_0}{\gamma}\right)\right],
\end{eqnarray*}
 where  $\sigma_2$ is defined such that $|\beta_{\sigma_2(1)}|\le \dots \le |\beta_{\sigma_2(k_0)}|$ and $A(k,t)$ is defined as in Theorem \ref{th_3}.\\
 Remark \ref{rem3} is also verified in the particular case where $(X_i)_{2\le i \le p}$ is an orthonormal family.
\end{remark}

\subsection{The high-dimensional case}
We will now discuss the non-ordered variable selection case in a high dimensional context, $p\ge n$. This section fits the two-step procedures previously introduced to high-dimensional analysis. \citet{Verzelen:2010} shows that when $k_0 log(p/k_0)/n\ge 1/2$, called the ultra-high-dimensional case, it is almost impossible to estimate the support of $\beta$. We will then consider that we are not in the ultra-high-dimensional case. 

The family $(X_i)_{1\le i\le p}$ is now a dependent family.
As said at the beginning of Section \ref{nonordered}, we assume that the decomposition of $\mu$ is unique, i.e. $\exists ! J\subset\{1,\dots,p\} /\mu=\sum_{j\in J}X_j\beta_j$. We still have $|J|=k_0$. The general procedure defined at the beginning of Section 3 remains the same: the first step orders the variables and the second step is a multiple testing of $\hat H_k$ defined by \eqref{Hk}.

Procedure `A' defined in Section \ref{sigyes} when the variance is known and procedure `B' defined in Section \ref{signot} when the variance is unknown are applicable with a minor modification. This modification concerns the definition of the subspaces $V_{(k)}$ and $S_{(k),(t)}$. Indeed, we have to take into account that $(X_i)_{1\le i\le p}$ is a linearly dependent family.

Let denote by $a_s$ the dimension of the linear span of $(X_{(1)},\dots,X_{(s)})$ for $1\le s \le p$, note that $a_p\le n$. 
We define, for all $1\le k< a_p-1$, $t^k_{max}=\lfloor log_2(a_p-k-1)\rfloor$ and $\mathcal{T}_k=\{0,\dots,t^k_{max} \}$. Denote $V_{k}=\text{span}(X_{(1)},\dots,X_{(s_k)})$ where $s_k$ is defined by $s_k=\text{inf}\{s: a_s=k\}$. 
As mentioned before, we test successively the hypotheses $\hat H_k: \{\mu \in V_{(k)}\}$ for all $1\le k<a_p-1$ until a null hypothesis is accepted. \\
In order to test $\hat H_k$, we consider the collection of alternatives $\{\mu \in V_{(k)}+S_{(k),(t)}\}$ for $t\in \mathcal{T}_k$ where \\$S_{(k),(t)}=\text{span}\left(\Pi_{V_{(k)}^\bot}(X_{(s_k+1)}),\dots,\Pi_{V_{(k)}^\bot}(X_{(s_k+q_{k,t})})\right)$\\ with $q_{k,t}$ defined by
$q_{k,t}=\text{inf}\{q:\text{dim}\left(Q_{(k),(q)}\right)=2^t\}$, where for all $q>0$,\\
$Q_{(k),(q)}=\text{span}\left(\Pi_{V_{(k)}^\bot}(X_{(s_k+1)}),\dots,\Pi_{V_{(k)}^\bot}(X_{(s_k+q)})\right)$.

With this modification, Lemma \ref{over12} and Lemma \ref{over22} still apply, thus Theorem \ref{th_2} of Section \ref{sigyes} and Theorem \ref{th_3} of Section \ref{signot} also apply. A simulation study is given in the next section showing that our procedure `B' performs well.

\section{Simulation study} 
\label{result}
\subsection{Presentation of the procedures}

In this section, we comment the results of the simulation study which are presented in tables \ref{orthon}-\ref{highdim}. Our aim was to compare the performances of our selection methods. The procedures presented in this paper are implemented in the R-package \textit{mht} which is available on CRAN (http://cran.r-project.org/).  Six methods were compared; the procedure described in Section \ref{ordered} for ordered variable selection, denoted ``proc-ordered" in the tables, the two-step procedure `B' described in Section \ref{nonordered}, either with ordered p-values denoted "procpval" or with the Bolasso order denoted ``procbol'', the FDR procedure described in \citet{bunea:2006}, the Lasso method and the Bolasso technique. For the purpose of comparison, we considered the design of the simulations of \citet{bunea:2006}. The comparison of the first method and the others is unfair and was not performed. The two kinds of method have to be separately compared.

The simulation was performed in several frameworks: in the common case when $(X_i)_{1\le i\le p}$ is a linearly independent family, in the orthonormal case, and in the high-dimensional case $(p\ge n)$. For the latter, the FDR procedure of \citet{bunea:2006} cannot be computed as p-values cannot be obtained with the least squares estimate  with all $p$ variables. In this case we compared an adjusted FDR (denoted FDR2); a p-value was calculated for each variable $X_i$ from the regression of Y onto that variable. As mentioned in the introduction, this is a natural extension of the FDR procedure in high-dimensional analysis and extended FDR is widely used in biology for QTL research and 
transcriptome analysis. 

When $(X_i)_{2\le i\le p}$ is not an orthonormal family, the calculation of $T_{k,\alpha}$ with \eqref{Tk} -for ordered variable selection- requires a lot of computational time, as a calculation of $V_k^\bot$ and $\{S_{k,t},t\in \mathcal{T}_k\}$ is needed for each $k$.  
Since a variable selection method is not only judged on its results but also on its fastness, useless calculations in our procedure had to be avoided. The Gram-Schmidt process was used to get an orthonormal family out of $(X_i)_{2\le i\le p}$. 
Thus the calculation of $(V_k^\bot)_{k\ge 0}$  was done once and for all.\\Decompose $\forall \ l>0$:
$X_{k+l}=\Pi_{V_k}(X{_{k+l}}) + \Pi_{V_k^{\bot}}(X_{k+l})$.
Note $(e_j)_{j=1,\dots,k}$ an orthonormal basis of $V_k$, then:\\
$\Pi_{V_k}(X{_{k+l}})=\sum_{j=1}^k <X_{k+l},e_j>e_j $
 and 
$\Pi_{V_k^{\bot}}(X_{k+l})=X_{k+l}-\sum_{j=1}^k <X_{k+l},e_j>e_j$.
The family $(X_1,\dots,X_p)$ was modified into $$\left(X_1,\dfrac{\Pi_{V_1^{\bot}}(X_{2})}{||\Pi_{V_1^{\bot}}(X_{2})    ||},\dfrac{\Pi_{V_2^{\bot}}(X_{3})}{||\Pi_{V_2^{\bot}}(X_{3})   ||},\dots,\dfrac{\Pi_{V_{p-1}^{\bot}}(X_{p})}{||\Pi_{V_{p-1}^{\bot}}(X_{p})    ||}\right) .$$
Denote that orthonormal family by $(\tilde{X}_1,\dots,\tilde{X}_p)$.
We decomposed $Y$ as:
\begin{equation}
\label{1bis}
Y=\underbrace{\tilde{X}_1\tilde{\beta}_1+\dots+\tilde{X}_k\tilde{\beta}_k}_{V_k}+\underbrace{\tilde{X}_{k+1}\tilde{\beta}_{k+1}+\dots+\tilde{X}_p\tilde{\beta}_p}_{\subset V_k^{\bot}}+\epsilon.
\end{equation}
Then $S_{k,t} =  \text{span}(\tilde{X}_{k+1},\dots,\tilde{ X}_{k+2^t})$ and so $\left|\left|\Pi_{S_{k,t}} Y\right|\right|_n^2=\tilde{\beta}_{k+1}^2+\dots+\tilde{\beta}_{k+2^t}^2$.
This technique avoided a lot of useless and redundant calculations.\\
The decomposition of Gram-Schmidt has also been used  in the non-ordered variable selection case with the two-step procedure `A' and `B' once the variables have been ordered.\\

When $(X_i)_{2\le i\le p}$ is an orthonormal family and the variance is unknown, we use another upper bound of the statistics $\tilde{U}_{D_{k,t},N_{k,t}}$ in our simulations than the one in Lemma \ref{over22}. Indeed, we can obtain an upper bound which does not depend on the family $(X_i)_{1\le i\le p}$ nor on the order on that family.\\
Let $I_1,\dots,I_{p}$ be $p$ i.i.d. standard Gaussian variables, and let $|I_{(1)}|\ge\dots\ge|I_{(p)}|$.\\ We define: $\forall k=0,\dots,p-1, \ \forall D=0,\dots,p-k-1$, $L_{k,D}=\sum_{j=k+D+1}^p I_{(j)}^2$.\\
Let $1\le k < p$ and $t\in \mathcal{T}_k$, we define \\$A_k=\{\{X_{(1)},\dots,X_{(k)}\}=\{X_j,j\in J\}\}$. For all $x>0$ we have

$$\mathbb{P}\left((\tilde{U}_{D_{k,t},N_{k,t}}>x)\cap A_k\right)\le \mathbb{P}\left(\dfrac{N_{k,t}}{D_{k,t}}\dfrac{Z_{D_{k,t},p-k}}{L_{k,D_{k,t}}+K_{n-p}}>x\right),$$
where $K_{n-p}$ is a chi-square variable with $n-p$ degrees of freedom and $Z_{d,D}$ is defined by \eqref{Zd}.
The simulations in the case of an orthonormal family and unknown variance were performed with this new upper bound. 

\subsection{Design of our simulation study}
Concerning the design of our simulations, we set $X_1$ to be the vector of $\mathbb{R}^n$ whose coordinates are all equal to $1/\sqrt{n}$ and we simulated $p-1$ independent vectors $X_j^* \sim \mathcal{N}_n(0,I_n)$, and set the predictors $X_j=X_j^*/||X_j^*||$, for $j=2,\dots,p$. The response variable $Y$ was computed via $Y=\beta_{i_1}X_{i_1}+\dots+\beta_{i_{k_0-1}}X_{i_{k_0-1}}+ \epsilon$, where $\epsilon$ is a vector of independent standard Gaussian variables, $J=\{1,i_{1},\dots,i_{k_0-1}\}\subset \{1,\dots,p\}$ and $\{\beta_{i_1},\dots,\beta_{i_{k_0-1}}\}\in\{\sqrt{n},6\}$.
We considered two instances of $k_0$ ($6$ or $11$). In each instance, samples of $n=100$ and $ 500$ in the low-dimensional case, and  $n=100$ in the high-dimen\-sional case have been simulated.
We let $p$ to vary with the sample size $n$.
In the orthonormal case,  we set the predictors $X_j,j=2,\dots,p$ as an orthonormal basis of $\text{span}(X_2^*,\dots,X_p^*)$ (principal component for example). 

In all tables, the first column gives an estimation of $\delta=P_{\mu}(A_{k_0}^c)$, calculated over $500$ replications. This estimation is not mentioned for the first procedure as we only considered ordered variable selection.
In low dimension, the parameter $m$ reflects whether the matrix $X$ is well-conditioned; $m=\text{max}_{1\le j\le p}\ m_{jj}$ where $(m_{ij})_{1\le i,j \le p}=(X^TX)^{-1}/n$, a low $m$ means that the matrix $X$ is well-conditioned.
The second column ``Truth" records the percentage of times the true model is selected; i.e. the pourcentage of time we actually found $\hat{J}=J$. The third column, labelled ``Inclusions", records the number of selected variables, averaged over 500 replications. ``Correct incl." records the number of relevant variables that are included in the selected model, averaged too. The MSE (Mean Squared Error) was calculated by average over all simulations:
$MSE=\sum_{i=1}^n (\hat{Y}_i-(X\beta_J)_i)     /n$, where $\hat{Y}=X\hat{\beta}$, where $\hat{\beta}$ is an estimation of $\beta$ with non zero values only on $\hat{J}$.

The FDR procedure described in \citet{bunea:2006} was set by choosing $q$ (user level) as $0.1$ and $0.05$. 

The $l^1$ penalty of the Lasso was tuned via 10-cross validation. Concerning the Bolasso technique, we choose $it=100$ bootstrap iterations; the frequency threshold and the pe\-nalty were also tuned via 10-cross validation. 
Both methods were always ended with a linear regression on the estimated set of indices of the relevant variables in order to minimize the bias of Lasso.

When the Bolasso technique was used to order the variable at the first step of procedure `B', we chose to stop the dichotomy algorithm (see Section \ref{nonordered}) as soon as $60$ variables were ordered. The objective was to spare calculation since it was uneasy to distinguish the remaining variables after the sixtieth position. The dichotomy algorithm assumes that when a variable reaches the frequency $1$, the frequency stays at $1$ when the penalty decreases. In practice this assumption might be wrong, the algorithm is then restarted.

Concerning the three procedures presented in this paper, the results are displayed for a level $\alpha \in \{0.1,0.05\}$. For ordered variable selection, $(X_i)_{1\le i\le p}$ was modified into $(X_J,X_{\{1,\dots,p\}\backslash J})$ and the collection $\{\alpha_{k,t}, t \in \mathcal{T}_k\}$ was chosen in accordance to the procedure P1, which required more computational time than P2, but which was far much powerful. For non-ordered variable selection, the collection $\{\alpha_{k,t}, t \in \mathcal{T}_k\}$ was chosen with the procedure P5 when $(X_i)_{2\le i\le p}$ was not an orthonormal family, since the variance was considered  unknown in the simulation.

\subsection{Comments on the results}
Table \ref{orthon} presents the results when the family $(X_i)_{1\le i\le p}$ is orthonormal, the non-orthonormal case with $p<n$ is reported in Table \ref{nonortho}. Lastly, Table \ref{highdim} shows the results of the methods in a high-dimensional context $(p\ge n)$; two alternatives were chosen for the number of variables, $p=300$ and $p=600$.

The procedure of multiple hypotheses testing developed for ordered variable selection in Section \ref{ordered} gave excellent results, even in the high-dimensional case where $p>n$, see Tables \ref{orthon}-\ref{highdim}. These results are not surprising because our choices of $\beta$ ensured that at each step the tests are powerful, so the probability of wrongly estimating $k_0$ was almost reduced to $\alpha$.

Concerning all the other tested methods, Table \ref{orthon} shows that the FDR procedure performed slightly better in the orthonormal case when $\beta_J$ was small. 
However, in the non-orthonormal case, our procedure `B' with the variables ordered thanks to the Bolasso technique gave the best results, especially compared to the FDR procedure which gave weak results (Table \ref{nonortho}).  

Table \ref{highdim} shows that the FDR2  was far from satisfactory. Indeed, nearly no true model were recovered in the $500$ simulations. 
 In fact, Table \ref{highdim} shows that our  ``procbol" procedure outperformed the others when $p>>n$. However, a combination of a small $\beta_J$  and a high number of variables induced a high $\hat{\delta}$ and consequently decreased the power of our ``procbol" method. 
Moreover, the results of the ``procbol" method become less satisfactory with an increase on the value of  $k_0$ because of the overestimation of the statistics in Lemma \ref{over22}.   \\

\begin{table}[h!]
\center
\label{orthon}
\footnotesize
\begin{tabular}{l| l l  | l l  | l l  | l   l|l | l }
Results &\multicolumn{2}{c|}{proc-ordered} & \multicolumn{2}{c|}{procpval} &\multicolumn{2}{c|}{procbol} & \multicolumn{2}{c|}{FDR} & Lasso& Bolasso\\
&$\alpha$=0.1 & $\alpha$=0.05 & $\alpha$=0.1 & $\alpha$=0.05 & $\alpha$=0.1 & $\alpha$=0.05  & $q$=0.1 & $q$=0.05 & \\ 
\hline
\multicolumn{5}{c}{$k_0=10,n=100,p=80,\beta=\sqrt{n},m=0.01$}\\
\hline
& & &\multicolumn{2}{c|}{$\hat{\delta}=0.00$}&\multicolumn{2}{c|}{$\hat{\delta}=0.00$}&\multicolumn{2}{c|}{$\hat{\delta}=0.00$}&\\
Truth 				&0.89	&0.95	&0.98	&0.99	&0.98	&0.99	&0.88	&0.95	&0.80	&0.75	\\
Inclusions 			&11.58	&11.22	&11.06	&11.01	&11.05	&11.02	&11.16	&11.07	&11.93	&11.80	\\
Correct incl.			&11.00	&11.00	&11.00	&10.99	&11.00	&11.00	&11.00	&11.00	&11.00	&11.00	\\
MSE		 		&0.11	&0.11	&0.11	&0.10	&0.10	&0.10	&0.11	&0.11	&0.15	&0.15	\\
\hline
\multicolumn{5}{c}{$k_0=5,n=100,p=80,\beta=6,m=0.01$}\\
\hline
& & &\multicolumn{2}{c|}{$\hat{\delta}=0.00$}&\multicolumn{2}{c|}{$\hat{\delta}=0.01$}&\multicolumn{2}{c|}{$\hat{\delta}=0.01$}&\\
Truth 				&0.89	&0.96	&0.75	&0.70	&0.80	&0.76	&0.81	&0.78	&0.73&0.72	\\
Inclusions 			&6.60	&6.19	&5.82	&5.68	&5.89	&5.78	&5.97	&5.80	&6.91&	6.82\\
Correct incl. 			&6.00	&6.00	&5.77	&5.67	&5.82	&5.74	&5.85	&5.74	&5.96& 5.99\\
MSE				&0.06	&0.06	&0.13	&0.16	&0.11	&0.14	&0.11	&0.15	&0.12&0.10\\
\end{tabular}
\caption{The orthonormal case. The first column gives an estimation of $\delta=P_{\mu}(A_{k_0}^c)$, calculated over $500$ replications. The second column ``Truth" records the pourcentage of time $\hat{J}=J$. ``Inclusions" records the number of selected variables and ``Correct incl." the number of selected variables that are relevant, averaged over 500 replications. The MSE  is calculated by average over all simulations:
$MSE=\sum_{i=1}^n (\hat{Y}_i-(X\beta_J)_i)     /n$, where $\hat{Y}=X\hat{\beta}$, where $\hat{\beta}$ is an estimation of $\beta$ with non zero values only on $\hat{J}$.
}\end{table}

\newpage
\begin{table}[h!]
\label{nonortho}
\center
\footnotesize
\begin{tabular}{l| l l  | l l  | l l  | l   l|l | l }
Results &\multicolumn{2}{c|}{proc-ordered} & \multicolumn{2}{c|}{procpval} &\multicolumn{2}{c|}{procbol} & \multicolumn{2}{c|}{FDR} & Lasso& Bolasso\\
&$\alpha$=0.1 & $\alpha$=0.05 & $\alpha$=0.1 & $\alpha$=0.05 & $\alpha$=0.1 & $\alpha$=0.05  & $q$=0.1 & $q$=0.05 & \\ 
\hline
\multicolumn{5}{c}{$k_0=10,n=100,p=80,\beta=\sqrt{n},m=0.102$}\\
\hline
& & &\multicolumn{2}{c|}{$\hat{\delta}=0.46$}&\multicolumn{2}{c|}{$\hat{\delta}=0.00$}&\multicolumn{2}{c|}{$\hat{\delta}=0.45$}&\\
Truth 				&0.92	&0.96	&0.54	&0.54	&0.94	&0.96	&0.13	&0.10	&0.29	&0.67	\\
Inclusions 			&11.33	&11.15	&13.06	&12.62	&11.08	&11.05	&8.55	&7.60	&13.18	&11.70	\\
Correct incl.		 	&11.00	&11.00	&10.92	&10.90	&11.00	&11.00	&8.34	&7.53	&11.00	&10.99	\\
MSE		 		&0.12	&0.11	&0.20	&0.22	&0.11	&0.11	&2.97	&3.72	&0.18	&0.14	\\
\hline
\multicolumn{5}{c}{$k_0=5,n=100,p=80,\beta=6,m=0.103$}\\
\hline
 & & &\multicolumn{2}{c|}{$\hat{\delta}=0.88$}&\multicolumn{2}{c|}{$\hat{\delta}=0.07$}&\multicolumn{2}{c|}{$\hat{\delta}=0.82$}&\\
Truth 				&0.91	&0.95	&0.11	&0.11	&0.86	&0.84	&0.00	&0.00	&0.27	&0.47	\\
Inclusions 			&6.37	&6.13	&7.30	&6.54	&6.00	&5.94	&1.98	&1.66	&8.22	&7.14	\\
Correct incl.		 	&6.00	&6.00	&5.05	&4.90	&5.91	&5.87	&1.86	&1.62	&5.94	&5.94	\\
MSE				&0.06	&0.06	&0.40	&0.44	&0.08	&0.09	&1.42	&1.45	&0.16	&0.13	\\
\hline
\multicolumn{5}{c}{$k_0=10,n=500,p=450,\beta=\sqrt{n},m=0.040$}\\
\hline
& & &\multicolumn{2}{c|}{$\hat{\delta}=0.02$}&\multicolumn{2}{c|}{$\hat{\delta}=0.00$}&\multicolumn{2}{c|}{$\hat{\delta}=0.01$}&\\
Truth		 		&0.91 	&0.95	&0.98	&0.98   	&0.94 	&0.96  	&0.84	&0.85	&0.88	&0.99	\\
Inclusions 			&12.09	&11.32	&11.05	&11.05	&11.07	&11.04	&11.12	&10.99	&11.26	&11.01	\\
Correct incl.		 	&11.00	&11.00	&11.00	&11.00	&11.00	&11.00	&10.94	&10.90	&11.00	&11.00	\\
MSE				&0.02	&0.02	&0.02	&0.02	&0.02	&0.02	&0.30	&0.31	&0.02	&0.02	\\
\hline
\multicolumn{5}{c}{$k_0=5,n=500,p=450,\beta=6,m=0.044$}\\
\hline
 & & &\multicolumn{2}{c|}{$\hat{\delta}=1.00$}&\multicolumn{2}{c|}{$\hat{\delta}=0.07$}&\multicolumn{2}{c|}{$\hat{\delta}=1.00$}&\\
Truth 				&0.89	&0.95	&0.00	&0.00	&0.86	&0.84	&0.00	&0.00	&0.68&	0.27	\\
Inclusions 			&8.35	&7.06	&3.19	&2.68	&5.95	&5.88	&1.09	&1.05	&6.37&	4.78\\
Correct incl.		 	&6.00	&6.00	&2.22	&2.16	&5.90	&5.85	&1.07	&1.05	&5.91&	4.78\\
MSE				&0.02	&0.01	&0.27	&0.28	&0.02	&0.02	&0.36	&0.36	&0.03	&0.09\\
\end{tabular} 
\caption{The non orthonormal case, $p< n$.  The first column gives an estimation of $\delta=P_{\mu}(A_{k_0}^c)$, calculated over $500$ replications. The second column ``Truth" records the pourcentage of time $\hat{J}=J$. ``Inclusions" records the number of selected variables and ``Correct incl." the number of selected variables that are relevant, averaged over 500 replications. The MSE  is calculated by average over all simulations:
$MSE=\sum_{i=1}^n (\hat{Y}_i-(X\beta_J)_i)     /n$, where $\hat{Y}=X\hat{\beta}$, where $\hat{\beta}$ is an estimation of $\beta$ with non zero values only on $\hat{J}$.}
\end{table}

\newpage
\begin{table}[h!]
\center
\label{highdim}
\footnotesize
\begin{tabular}{l| l l  | l l  | l l  | l   l|l | l }
Results &\multicolumn{2}{c|}{proc-ordered} & \multicolumn{2}{c|}{procpval} &\multicolumn{2}{c|}{procbol} & \multicolumn{2}{c|}{FDR2} & Lasso& Bolasso\\
&$\alpha$=0.1 & $\alpha$=0.05 & $\alpha$=0.1 & $\alpha$=0.05 & $\alpha$=0.1 & $\alpha$=0.05  & $q$=0.1 & $q$=0.05 & \\ 
\hline
\multicolumn{5}{c}{$k_0=10,n=100,p=300,\beta=\sqrt{n}$}\\
\hline
& & &\multicolumn{2}{c|}{$\hat{\delta}=1.00$}&\multicolumn{2}{c|}{$\hat{\delta}=0.00$}&\multicolumn{2}{c|}{$\hat{\delta}=1.00$}&\\
Truth 				&0.91	&0.96	&0.00	&0.00	&0.99	&0.99	&0.00	&0.00	&0.60	&0.78	\\
Inclusions 			&11.53	&11.14	&9.92	&9.68	&11.01	&11.01	&5.17	&4.38	&12.05	&11.46	\\
Correct incl.		 	&11.00	&11.00	&9.35	&9.36	&11.00	&11.00	&5.17	&4.38	&11.00	&11.00\\
MSE		 		&0.11	&0.10	&1.56	&1.63	&0.10	&0.10	&5.21	&6.04	&0.15	&0.13\\
\hline
\multicolumn{5}{c}{$k_0=5,n=100,p=300,\beta=6$}\\
\hline
& & &\multicolumn{2}{c|}{$\hat{\delta}=0.65$}&\multicolumn{2}{c|}{$\hat{\delta}=0.09$}&\multicolumn{2}{c|}{$\hat{\delta}=0.60$}&\\
Truth 				&0.93	&0.96	&0.33	&0.33	&0.79	&0.74	&0.03	&0.01	&0.38	&0.56	\\
Inclusions 			&6.44	&6.16	&5.62	&5.50	&5.88	&5.78	&4.22	&3.74	&8.57	&7.32	\\
Correct incl.		 	&6.00	&6.00	&5.32	&5.24	&5.82	&5.74	&4.15	&3.71	&5.92	&5.90	\\
MSE		 		&0.06	&0.05	&0.27	&0.29	&0.11	&0.14	&0.66	&0.79	&0.18	&0.15	\\
\hline
\multicolumn{5}{c}{$k_0=10,n=100,p=600,\beta=\sqrt{n}$}\\
\hline
& & &\multicolumn{2}{c|}{$\hat{\delta}=1.00$}&\multicolumn{2}{c|}{$\hat{\delta}=0.17$}&\multicolumn{2}{c|}{$\hat{\delta}=1.00$}&\\
Truth 				&0.89	&0.95	&0.00	&0.00	&0.83	&0.83	&0.00	&0.00	&0.00	&0.25	\\
Inclusions 			&11.66	&11.21	&5.88	&5.36	&11.30	&11.20	&3.33	&3.02	&17.97	&13.24	\\
Correct incl.		 	&11.00	&11.00	&5.68	&5.23	&10.99	&10.99	&3.33	&3.02	&10.99	&10.99	\\
MSE		 		&0.12	&0.11	&4.11	&4.56	&0.11	&0.11	&6.34	&6.69	&0.31	&0.20	\\
\hline
\multicolumn{5}{c}{$k_0=5,n=100,p=600,\beta=6$}\\
\hline
& & &\multicolumn{2}{c|}{$\hat{\delta}=0.95$}&\multicolumn{2}{c|}{$\hat{\delta}=0.30$}&\multicolumn{2}{c|}{$\hat{\delta}=0.92$}&\\
Truth 				&0.91	&0.96	&0.05	&0.05	&0.62	&0.56	&0.00	&0.00	&0.11	&0.26	\\
Inclusions 			&6.43	&6.12	&4.36	&4.22	&5.62	&5.48	&2.48	&2.18	&11.52	&8.49	\\
Correct incl.		 	&6.00	&6.00	&4.14	&4.04	&5.50	&5.39	&2.46	&2.17	&5.59	&5.65	\\
MSE		 		&0.06	&0.06	&0.59	&0.62	&0.22	&0.25	&1.10	&1.22	&0.37	&0.30	\\
\end{tabular}
\caption{The high-dimensional case, $p\ge n$. The first column gives an estimation of $\delta=P_{\mu}(A_{k_0}^c)$, calculated over $500$ replications. The second column ``Truth" records the pourcentage of time $\hat{J}=J$. ``Inclusions" records the number of selected variables and ``Correct incl." the number of selected variables that are relevant, averaged over 500 replications. The MSE  is calculated by average over all simulations:
$MSE=\sum_{i=1}^n (\hat{Y}_i-(X\beta_J)_i)     /n$, where $\hat{Y}=X\hat{\beta}$, where $\hat{\beta}$ is an estimation of $\beta$ with non zero values only on $\hat{J}$.
}\end{table}

\section{Conclusion }

This paper tackled the problem of recovering the set of relevant variables $J$ in a sparse linear model, especially when the number of variables $p$ was higher than the sample size $n$. We proposed new methods based on hypotheses testing to estimate $J$.
The procedures are proved to be powerful under some conditions on the signal and the theorems are non asymptotic.
The simulations showed that these new procedures outperformed all the other tested methods, especially in the high-dimensional case, which was the aim of this study. 

\bibliography{test}

\begin{thebibliography}{}

\bibitem[Bach, 2009]{bach:2009}
Bach, F. (2009).
\newblock Model-consistent sparse estimation through the bootstrap.

\bibitem[Baraud et~al., 2003]{Baraud:2003}
Baraud, Y., Huet, S., and Laurent, B. (2003).
\newblock Adaptative test of linear hypotheses by model selection.
\newblock {\em Ann. Statist.}, 31(1):225--251.

\bibitem[Benjamini and Hochberg, 1995]{Benjamini:1995}
Benjamini, Y. and Hochberg, Y. (1995).
\newblock Controlling the false discovery rate: a practical and powerful
  approach to multiple hypothesis testing.
\newblock {\em J. R. Stat. Soc.}, B 57, 289-300.

\bibitem[Bickel et~al., 2009]{Bickel:2009}
Bickel, P.~J., Ritov, Y., and Tsybakov, A.~B. (2009).
\newblock Simultaneous analysis of lasso and dantzig selector.
\newblock {\em Ann. Statist.}, 37(4):1705--1732.

\bibitem[Bunea et~al., 2007]{Bunea:2007}
Bunea, F., Tsybakov, A., and Wegkamp, M. (2007).
\newblock Sparsity oracle inequalities for the lasso.
\newblock {\em Electron. J. Statist.}, 1:169--194.

\bibitem[Bunea et~al., 2006]{bunea:2006}
Bunea, F., Wegkamp, M., and Auguste, A. (2006).
\newblock Consistent variable selection in high dimensional regression via
  multiple testing.
\newblock {\em Statist. Plann. Inference}, 136:4349--4363.

\bibitem[Candes and Tao, 2007]{Candes:2007}
Candes, E. and Tao, T. (2007).
\newblock The dantzig selector: Statistical estimation when p is much larger
  than n.
\newblock {\em Ann. Statist.}, 35(6):2313--2351.

\bibitem[Chesneau and Hebiri, 2008]{Chesneau:2008}
Chesneau, C. and Hebiri, M. (2008).
\newblock Some theoretical results on the grouped variables lasso.
\newblock {\em Math. Methods Statist.}, 17(4):317--326.

\bibitem[Huang et~al., 2008]{Huang:2008}
Huang, J., Ma, S., and Zhang, C.-H. (2008).
\newblock Adaptative lasso for sparse high-dimensional regression models.
\newblock {\em Stat. Sin.}, 18(4):1603--1618.

\bibitem[Laurent and Massart, 2000]{Laurent:2000}
Laurent, B. and Massart, P. (2000).
\newblock Adaptive estimation of a quadratic functional by model selection.
\newblock {\em Ann. Statist.}, 28(5):1302--1338.

\bibitem[L\^e~Cao et~al., 2008]{LeCao:2008}
L\^e~Cao, K.~A., Rossouw, D., Robert-Grani\'e, C., and Besse, P. (2008).
\newblock A sparse pls for variable selection when integrating omics data.
\newblock {\em Statistical applications in genetics and molecular biology},
  7:Article 35.

\bibitem[Meinshausen and B\"uhlmann, 2006]{Meinshausen:2006}
Meinshausen, N. and B\"uhlmann, P. (2006).
\newblock High-dimensional graphs and variable selection with the lasso.
\newblock {\em Ann. Statist.}, 34(3):1436--1462.

\bibitem[Tenenhaus, 1998]{Tenenhaus:1998}
Tenenhaus, M. (1998).
\newblock {\em La r\'egression PLS: th\'eorie et pratique}.
\newblock Editions Technip.

\bibitem[Tibshirani, 1996]{Tibshirani:1996}
Tibshirani, R. (1996).
\newblock Regression shrinkage and selection via the lasso.
\newblock {\em J. R. Stat. Soc.}, B 58(1):267--288.

\bibitem[Verzelen, 2012]{Verzelen:2010}
Verzelen, N. (2012).
\newblock Minimax risks for sparse regressions: Ultra-high-dimensional
  phenomenons.
\newblock {\em Electron. J. Statist.}, 6(1):38--90.

\bibitem[Zhao and Yu, 2006]{Zhao:2006}
Zhao, P. and Yu, B. (2006).
\newblock On model selection consistency of lasso.
\newblock {\em J. Mach. Learn. Res.}, 7:2541--2563.

\end{thebibliography}

\appendix
\section{Proofs}
\label{proof}
\begin{proof}[\textbf{Proof of Theorem \ref{th_1}}]

Let $k\le k_0-1$ and assume that $(R_k)$ holds. According to \citet{Baraud:2003}, the power of the test $H_{k}$, $\mathbb{P}_{\mu }(T_{k,\alpha}>0)$, is greater than $ 1- \gamma/k_0$. This is equivalent to \\$\mathbb{P}_{\mu }(H_k \ \text{is accepted}) \le \gamma/k_0$.\\
Moreover, for all $ k\ge k_0$, $\mathbb{P}_{\mu}(T_{k,\alpha}>0)\le \alpha$, since $\alpha$ is the level of the test $H_{k}$.
Then we have:
\begin{eqnarray*}
\mathbb{P}_{\mu }(\hat{k}> k_0) &\le& \mathbb{P}_{\mu }(H_{k_0} \text{ is rejected})
=\mathbb{P}_{\mu }(T_{k_0,\alpha}>0)\\
 &\le& \alpha
 \end{eqnarray*}
 and
\begin{eqnarray*}
\mathbb{P}_{\mu }(\hat{k}< k_0) & \le &  \sum_{j=0}^{k_0-1} \mathbb{P}_{\mu }(H_j \ \text{is accepted})\\
 & \le & k_0  \gamma/k_0.
\end{eqnarray*}

Hence we obtain
\begin{equation*}
\mathbb{P}_{\mu }(\hat{k}\neq k_0) \le \mathbb{P}_{\mu }(\hat{k}< k_0) +\mathbb{P}_{\mu }(\hat{k}> k_0) \le   \gamma + \alpha,
\end{equation*}
which concludes the proof of \eqref{non}.\\
\end{proof}

 \begin{proof}[\textbf{Proof of Lemma \ref{over12}}]
Let $x>0$. By definition of  $U_{k,t}$, we have\\
$\mathbb{P}\left((U_{k,t}>x)\cap A_k\right)$
\begin{eqnarray*}
\hspace{0.5cm}&=&\mathbb{P}\left(\left(\dfrac{||\Pi_{S_{(k),(t)}} Y||_n^2}{\sigma^2}>x\right)\cap A_k\right)\\
&=&\mathbb{P}\left(\left(\dfrac{||\Pi_{S_{(k),(t)}} \mu||_n^2+||\Pi_{S_{(k),(t)}} \epsilon||_n^2}{\sigma^2}>x\right)\cap A_k\right).
\end{eqnarray*}
Since $A_k=\{\{X_{(1)},\dots,X_{(k)}\}=\{X_j,j\in J\}\}$, \\
$\mathbb{P}\left(\left(\dfrac{||\Pi_{S_{(k),(t)}} \mu||_n^2+||\Pi_{S_{(k),(t)}} \epsilon||_n^2}{\sigma^2}>x\right)\cap A_k\right)$
\begin{eqnarray*}
\hspace{0.5cm}&=&\mathbb{P}\left(\left(\dfrac{||\Pi_{S_{(k),(t)}} \epsilon||_n^2}{\sigma^2}>x\right)\cap A_k\right)\\
&\le&\mathbb{P}\left(\dfrac{||\Pi_{S_{(k),(t)}} \epsilon||_n^2}{\sigma^2}>x\right).
\end{eqnarray*}
And by construction of $U_{k,t}^1$,
$$\mathbb{P}\left(\dfrac{||\Pi_{S_{(k),(t)}} \epsilon||_n^2}{\sigma^2}>x\right)\le\mathbb{P}\left(U_{k,t}^1>x\right).$$
Thus $$\mathbb{P}\left((U_{k,t}>x)\cap A_k\right)\le \mathbb{P}\left(U_{k,t}^1>x\right).$$
  \end{proof}

\begin{proof}[\textbf{Proof of Theorem \ref{th_2}}]
Let $k<k_0$.\\
We use the identity $\forall (a,b)\in \mathbb{R}^2, (a+b)^2\ge \dfrac{1}{2}a^2-b^2$.\\ On the event $A_{k_0}$:\\
$\forall t\in I= \{0,...,log_2(k_0-k)\}$:
\begin{eqnarray*}
||\Pi_{S_{(k),(t)}}Y ||_n^2&=&||\Pi_{S_{(k),(t)}}(\mu+\epsilon) ||_n^2\\
&\ge& \dfrac{1}{2}|| \Pi_{S_{(k),(t)}}\mu||_n^2- || \Pi_{S_{(k),(t)}}\epsilon||_n^2\\
&\ge& \dfrac{1}{2} \text{inf}\{|| \Pi_{S}\mu||_n^2,S\in B_{2^t}\}- || \Pi_{S_{(k),(t)}}\epsilon||_n^2
\end{eqnarray*}
where $B_{2^t}=\{\text{span}(X_I),I\subset J, |I|=2^t\}$.
Hence: \\

$\mathbb{P}\left(\forall t\in I, \dfrac{1}{\sigma^2}||\Pi_{S_{(k),(t)}}Y ||_n^2 \le  \overline{U^1_{k,t}}^{-1}(\alpha_{k,t}) \ \cap A_{k_0}\right) $
\begin{eqnarray*}
&=&\mathbb{P}\left( \forall t\in I,\dfrac{1}{\sigma^2}||\Pi_{S_{(k),(t)}}(\mu+\epsilon) ||_n^2 \le \overline{U^1_{k,t}}^{-1}(\alpha_{k,t}) \ \cap  A_{k_0}\right)\\
&\le& \mathbb{P}\left(\forall t\in I,\dfrac{1}{2\sigma^2} \text{inf}\{|| \Pi_{S}\mu||_n^2,S\in B_{2^t}\}- \dfrac{1}{\sigma^2}|| \Pi_{S_{(k),(t)}}\epsilon||_n^2\le\overline{U^1_{k,t}}^{-1}(\alpha_{k,t})\ \cap  A_{k_0}\right).
\end{eqnarray*}

We have on the event $A_{k_0}$ and for $k+2^t\le k_0$ that\\ $|| \Pi_{S_{(k),(t)}}\epsilon||_n^2\le \text{sup}\{|| \Pi_{S}\epsilon||_n^2,S\in B_{2^t}\}$. Moreover, for $S\in B_{2^t}$, $|| \Pi_{S}\epsilon||_n^2\sim \dfrac{\sigma^2}{n}\chi^2_{2^t}$. Note that $|B_{2^t}|=\begin{pmatrix}k_0\\2^t\end{pmatrix}$. \\
Denote $Z_t=\dfrac{|| \Pi_{S_{(k),(t)}}\epsilon||_n^2}{\sigma^2}$ and $\bar{Z}_t(u)$ the probability for the statistic $Z_t$ to be larger than u. We denote $\bar{\chi}_d(u)$ the probability for a chi-square with $d$ degrees of freedom to be larger than $u$.
We have an upper bound of the $(1-u)$-quantile of the statistic $Z_t$: $\bar{Z}_t^{-1}(u)\le \bar{\chi}_{2^t}^{-1}(u/|B_{2^t}|)/n$.    Indeed:\\
\begin{eqnarray*}
\mathbb{P}\left(Z_t>\dfrac{\bar{\chi}_{2^t}^{-1}(u/|B_{2^t}|)}{n}\right) &\le &\mathbb{P}\left( \text{sup}\{\dfrac{ || \Pi_{S}\epsilon||_n^2}{\sigma^2},S\in B_{2^t}\}>\dfrac{\bar{\chi}_{2^t}^{-1}(u/|B_{2^t}|)}{n}\right)\\
&\le& \sum_{S\in B_{2^t}} \mathbb{P}\left(|| \Pi_{S}\epsilon||_n^2>\dfrac{\sigma^2}{n}\bar{\chi}_{2^t}^{-1}(u/|B_{2^t}|)\right)\\
&\le&|B_{2^t}|\dfrac{u}{|B_{2^t}|}\le u.
\end{eqnarray*}
Therefore, the following condition\\
$(cond_k):\exists t\in I,$
\begin{equation*}
\dfrac{1}{2\sigma^2} \text{inf}\{|| \Pi_{S}\mu||_n^2,S\in B_{2^t}\}\ge\dfrac{1}{n}\bar{\chi}_{2^t}^{-1}\left(\dfrac{\gamma/k_0}{|B_{2^t}|}\right)+\overline{U^1_{k,t}}^{-1}(\alpha_{k,t})
\end{equation*}
implies that:
\begin{equation}
\label{imply}
\mathbb{P}\left[\forall t\in I,\dfrac{1}{2\sigma^2} \text{inf}\{|| \Pi_{S}\mu||_n^2,S\in B_{2^t}\}- \dfrac{|| \Pi_{S_{(k),(t)}}\epsilon||_n^2}{\sigma^2}\le\overline{U^1_{k,t}}^{-1}(\alpha_{k,t})\ \cap  A_{k_0}\right]\le \gamma/k_0.
\end{equation}
Let us denote $\forall 0<d, $ 
 \begin{equation}
 \label{Gkd}
 G_{k,d}=\{\text{span}(X_I), I\subset \{1,..,p\}\backslash\{(1),...,(k)\},|I|=d \}.
 \end{equation}\\
 Note that $|G_{k,d}|=\begin{pmatrix}p-k\\d\end{pmatrix}$.\\ Then $U_{k,t}^1\le \text{sup}\{|| \Pi_{S}\epsilon||_n^2,S\in G_{k,2^t}\}$. This inequality leads us to an upper bound of the (1-u)-quantile of $U^1_{k,t}$:\\ $\overline{U_{k,t}^1}^{-1}(u)\le \bar{\chi}_{2^t}^{-1}(u/|G_{k,2^t}|)/n$.\\\\
Using $\overline{U_{k,t}^1}^{-1}(u)\le \bar{\chi}_{2^t}^{-1}(u/|G_{k,2^t}|)/n$ in the condition $(cond_k)$, we obtain the following condition which still implies \eqref{imply}:\\
 $$\exists t\in I,
 \dfrac{1}{2\sigma^2} \text{inf}\{|| \Pi_{S}\mu||_n^2,S\in B_{2^t}\}\ge\dfrac{1}{n}\left[\bar{\chi}_{2^t}^{-1}\left(\dfrac{\gamma/k_0}{|B_{2^t}|}\right)+\bar{\chi}_{2^t}^{-1}\left(\dfrac{\alpha_{k,t}}{|G_{k,2^t}|}\right)\right].$$
Moreover, \citet{Laurent:2000} showed that for $K\sim\chi^2_d$:
\begin{equation}
\label{chi2}
\forall x>0, \mathbb{P}\left({K}\ge d+2\sqrt{dx}+2x  \right)\le e^{-x}.
\end{equation}
Then for $d=2^t$ and  $x_u=log\left(\dfrac{|B_{2^t}|}{\gamma/k_0}\right)$ we have $\bar{\chi}_{2^t}^{-1}\left(\dfrac{\gamma/k_0}{|B_{2^t}|}\right)\le 2^t+2\sqrt{2^tx_u}+2x_u$. Since $\begin{pmatrix}D\\d\end{pmatrix}\le \left(\dfrac{eD}{d}\right)^d$, $|B_{2^t}|\le \left(\dfrac{ek_0}{2^t}\right)^{2^t}$, thus $x_u=2^t log\left(\dfrac{ek_0}{2^t}\right)+log\left(\dfrac{k_0}{\gamma}\right)$.\\
Using $\sqrt{u+v}\le \sqrt{u}+\sqrt{v}$ for all $u>0,v>0$ and $\sqrt{u}\le u$ for all $u\ge 1$,
we obtain: \\
\begin{eqnarray*}
\bar{\chi}_{2^t}^{-1}\left(\dfrac{\gamma/k_0}{|B_{2^t}|}\right)&\le& 2^t\left[1+2\sqrt{log\left(\dfrac{ek_0}{2^t}\right)}+2log\left(\dfrac{ek_0}{2^t}\right)\right]\\&&+2\left[\sqrt{2^t log(k_0/\gamma)}+log(k_0/\gamma)\right]\\
&\le& 2^t\left[5 +4log\left(\dfrac{k_0}{2^t}\right)\right]+2\left[\sqrt{2^t log(k_0/\gamma)}+log(k_0/\gamma)\right].
\end{eqnarray*}
For $d=2^t$ and $x_u=log(|G_{k,2^t}|/\alpha_{k,t})$, we obtain: \\
$\bar{\chi}_{2^t}^{-1}\left(\alpha_{k,t}/|G_{k,2^t}|\right)$
\begin{eqnarray*}
&\le& 2^t\left[1+2\sqrt{log\left(\dfrac{e(p-k)}{2^t}\right)}+2log\left(\dfrac{e(p-k)}{2^t}\right)\right]\\
&&+2\left[\sqrt{2^t log(1/\alpha_{k,t})}+log(1/\alpha_{k,t})\right]\\
&\le&2^t\left[5 +4log\left(\dfrac{p-k}{2^t}\right)\right]+2\left[\sqrt{2^t log(1/\alpha_{k,t})}+log(1/\alpha_{k,t})\right].
\end{eqnarray*}
We also have an upper bound of $1/\alpha_{k,t}, \forall t\in \mathcal{T}_k$. Indeed , the construction of $\{\alpha_{k,t}, t\in \mathcal{T}_k\}$ with the procedure P3 gives $\mathbb{P}\left(\exists t\in \mathcal{T}_k, U_{k,t}^1>\overline{U_{k,t}^1}^{-1}(\alpha_{k,t})\right)=\alpha$.
Thus $\forall t\in \mathcal{T}_k, \alpha_{k,t}\ge \alpha/|\mathcal{T}_k|$, since $\mathbb{P}\left(\exists t\in \mathcal{T}_k, U_{k,t}^1>\overline{U_{k,t}^1}^{-1}(\alpha/|\mathcal{T}_k|)\right)\le\alpha$.\\
Hence we obtain:
\begin{eqnarray*}
\bar{\chi}_{2^t}^{-1}\left(\alpha_{k,t}/|G_{k,2^t}|\right)&\le& 2^t\left[5 +4log\left(\dfrac{p-k}{2^t}\right)\right]\\&&+2\left[\sqrt{2^t log\left(\dfrac{|\mathcal{T}_k|}{\alpha}\right)}+ log\left(\dfrac{|\mathcal{T}_k|}{\alpha}\right)\right].
\end{eqnarray*}
Using the inequality $a\sqrt{u}+b\sqrt{v}\le \sqrt{a^2+b^2}\sqrt{u+v}$ which holds for any positive numbers $a, b, u, v$, we finally get the condition $(R_{2,k})$ which implies \eqref{imply}:\\
$(R_{2,k}):\ \exists t\in I$ such that\\
\begin{eqnarray*}
\dfrac{1}{2\sigma^2} \text{inf}\{|| \Pi_{S}\mu||_n^2,S\in B_{2^t}\}&\ge&\dfrac{2^t}{n} \left[10 +4log\left(\dfrac{(p-k)k_0}{2^{2t}}\right)\right]\\
&+&\dfrac{2}{n}\left[\sqrt{2^{t+1} log\left(\dfrac{k_0|\mathcal{T}_k|}{\gamma\alpha}\right) }+ log\left(\dfrac{k_0|\mathcal{T}_k|}{\gamma\alpha}\right) \right].
\end{eqnarray*}
This leads to
$$\mathbb{P}\left(\forall t\in I, \dfrac{1}{\sigma^2}||\Pi_{S_{(k),(t)}}Y ||_n^2\le \overline{U^1_{k,t}}^{-1}(\alpha_{k,t})\ \cap  A_{k_0}\right) \le \gamma/k_0.$$
Hence $$ \mathbb{P}\left(\forall t\in I, U_{k,t}\le\overline{U^1_{k,t}}^{-1}(\alpha_{k,t})\ \cap  A_{k_0}\right) \le \gamma/k_0.$$
Then, $\forall k<k_0, \mathbb{P}\left(\mathring{k}_{A}=k \ \cap \ A_{k_0}\right)\le\gamma/k_0$.

We can calculate $\mathbb{P}_{\mu}(\hat{J}\neq J )$:
\begin{eqnarray*}
\mathbb{P}_{\mu}(\hat{J}\neq J )& \le&\mathbb{P}_{\mu}(\hat{J}\neq J \ \cap  A_{k_0})+\mathbb{P}( A_{k_0}^c) \\
&\le &\left(\sum_{j=0}^{k_0-1} \mathbb{P}_{\mu}(\mathring{k}_{A}=j \ \cap  A_{k_0})+\mathbb{P}_{\mu}(\mathring{k}_{A}>k_0 \ \cap  A_{k_0})\right)+\mathbb{P}( A_{k_0}^c) \\
&\le& k_0  \gamma/k_0 + \alpha + \delta.
\end{eqnarray*}

And then \eqref{non1} is proved.
\end{proof}
 
 \begin{proof}[\textbf{Proof of Lemma \ref{over1}}]
 Under $\hat{H}_k$ and on the event $A_{k}$ :
 \begin{eqnarray*}
{U}_{k,t}&=&||\Pi_{S_{(k),(t)}} Y||_n^2/\sigma^2=||\Pi_{S_{(k),(t)}}(\mu+\epsilon) ||_n^2/\sigma^2\\
 &=&||\Pi_{S_{(k),(t)}}\epsilon ||_n^2/\sigma^2.
 \end{eqnarray*}
The family $(X_i)_i$ is orthonormal, thus:\\
$U_{k,t}=\sum_{j=k+1}^{k+2^t}<\epsilon,{X}_{(j)}>_n^2/\sigma^2$.\\
As $\epsilon\sim\mathcal{N}_n(0,\sigma^2I_n)$, we have for all $1\le j\le p$,\\ $<\epsilon,{X}_{j}>\sim \mathcal{N}(0,\sigma^2)$
and the variables $<\epsilon,X_j>,j=1,...,p$ are i.i.d..
Thus 
$\{<\epsilon,X_{(j)}>,j>k\}=\{<\epsilon,X_m>,m\notin J\}=\{\sigma W_1,...,\sigma  W_{p-k}\}$.

So
$\sum_{j=k+1}^{k+2^t}<\epsilon,{X}_{(j)}>_n^2/\sigma^2\le \sum_{j=1}^{2^t} W_{(j)}^2/n = Z_{k,D_{k,t}}/n$.\\
\end{proof}

\begin{proof}[\textbf{Proof of Corollary \ref{cor_2}}]
Let $k<k_0$.\\
$\sigma_2$ is defined such that $|\beta_{\sigma_2(1)}|\le ... \le |\beta_{\sigma_2(k_0)}|$, 
note $\epsilon_{(j+1)}=|| \Pi_{S_{(j),0}}\epsilon||$, $\forall j \in \{k+1,...,k+2^t\} \text{ with } k+2^t\le k_0$.\\
Similarly as in the proof of Theorem \ref{th_2}, using that \\$\text{inf}\{|| \Pi_{S}\mu||_n^2,S\in B_{2^t}\}=\sum_{j=1}^{2^t} \beta^2_{\sigma_2(j)},$ we get that:\\
$\mathbb{P}\left(\forall t\in I, \dfrac{1}{\sigma^2}||\Pi_{S_{(k),(t)}}Y ||_n^2 \le  \dfrac{\bar{Z}_{D_{k,t},p-k}^{-1}(\alpha_{k,t})}{n} \ \cap A_{k_0}\right) $
\begin{eqnarray*}
\hspace{1cm}&\le& \mathbb{P}\left(\forall t\in I,\dfrac{1}{2\sigma^2}\sum_{j=1}^{2^t} \beta^2_{\sigma_2(j)}-\dfrac{1}{n\sigma^2}\sum_{j=k}^{k+2^t-1} \epsilon_{(j+1)}^2\le\dfrac{\bar{Z}_{D_{k,t},p-k}^{-1}(\alpha_{k,t})}{n}\ \cap  A_{k_0}\right).
\end{eqnarray*}
On the event $A_{k_0}$, 
$\{<\epsilon, X_{(j+1)}>,k\le j\le k+2^t-1\}\subset\{<\epsilon,X_j>,j\in J\}$, which implies that we have an stochastic upper bound: $\sum_{j=k}^{k+2^t-1} \epsilon_{(j+1)}^2\le \sigma^2Z_{2^t,k_0}$. \\ 
Hence the following condition \\
\begin{equation}
\label{cond3}
\exists \ t \le log_2(k_0-k) /\dfrac{1}{2\sigma^2}\sum_{j=1}^{2^t} \beta_{\sigma_2(j)}^2 \ge   \dfrac{1}{n}\left[\bar{Z}_{D_{k,t},p-k}^{-1}(\alpha_{k,t})+\bar{Z}_{D_{k,t},k_0}^{-1}(\gamma/k_0) \right]
\end{equation}
implies that$$\mathbb{P}\left(\forall t\in I,\dfrac{1}{2\sigma^2}\sum_{j=1}^{2^t} \beta^2_{\sigma_2(j)}-\dfrac{1}{n\sigma^2}\sum_{j=k}^{k+2^t-1} \epsilon_{(j+1)}^2\le\dfrac{\bar{Z}_{D_{k,t},p-k}^{-1}(\alpha_{k,t})}{n}\ \cap  A_{k_0}\right)$$
$$\le \gamma/k_0.$$
This leads to
\begin{equation}
\label{imply2}
\mathbb{P}\left(\forall t\in I, \dfrac{1}{\sigma^2}||\Pi_{S_{(k),(t)}}Y ||_n^2\le \bar{Z}_{D_{k,t},p-k}^{-1}(\alpha_{k,t})\ \cap  A_{k_0}\right) \le \gamma/k_0.
\end{equation}
Let $0<u<1$, $0<D$ and $d<D$. 
In the following, we study the behavior of the $(1-u)$ quantile of the statistic $Z_{d,D}$ in order to obtain a more explicit condition than \eqref{cond3}.\\
 Let define $V_{d,D}=\{I \subset\{1,...,D\} / |I|=d\}$. Note that $|V_{d,D}|= \begin{pmatrix}D\\d\end{pmatrix}$.
Let recall that $Z_{d,D}$ is defined by \eqref{Zd} as $Z_{d,D}=\sum_{j=1}^d W_{(j)}^2$ where $W_1,...,W_{D}$ are $D$ i.i.d. standard Gaussian variables ordered as $|W_{(1)}|\ge...\ge|W_{(D)}|$.\\
We have that: $Z_{d,D}\le \text{sup}\{\sum_{i\in I} W_{i}^2,  I \in V_{d,D}\}$. Note that for $ I \in V_{d,D}, \sum_{i\in I} W_{i}^2,\sim \chi^2_d$.\\
We obtain that the $(1-u)$-quantile of $Z_{d,D}$ is lower than $\bar{\chi_d}^{-1}\left(u/|V_{d,D}|\right)$:\\
\begin{eqnarray*}
\mathbb{P}\left(Z_{d,D}>\bar{\chi_d}^{-1}\left(u/|V_{d,D}|\right)\right)&\le &\mathbb{P}\left(\text{sup}\{\sum_{i\in I} W_{i}^2, \forall I \in V_{d,D}\}>\bar{\chi_d}^{-1}\left(u/|V_{d,D}|\right)\right)\\
&\le& \sum_{I \in V_{d,D}} \mathbb{P}\left(\sum_{i\in I} W_{i}^2>\bar{\chi_d}^{-1}\left(u/|V_{d,D}|\right)\right)\\
&\le& |V_{d,D}|\dfrac{u}{|V_{d,D}|}\le u.
\end{eqnarray*}
Using the expression of the upper bound of $\bar{\chi}_d^{-1}(u)$ from the proof of Theorem \ref{th_2}, we get the condition $(R_{2bis,k})$ from an upper bound of the right part in the condition \eqref{cond3}. The end of the proof is the same as in the proof of Theorem \ref{th_2}.
\end{proof}

\begin{proof}[\textbf{Proof of Lemma \ref{over22}}]
Let $x>0$. By definition of  $U_{k,t}$, we have\\
\begin{eqnarray*}
\mathbb{P}\left((\tilde{U}_{D_{k,t},N_{k,t}}>x)\cap A_k\right)&=&\mathbb{P}\left(\left(\dfrac{N_{k,t} ||\Pi_{S_{(k),(t)}} Y||_n^2}{D_{k,t} ||Y-\Pi_{V_{(k),(t)}}Y||_n^2}>x\right)\cap A_k\right)\\
&=&\mathbb{P}\left(\left(\dfrac{N_{k,t}||\Pi_{S_{(k),(t)}} \mu||_n^2+N_{k,t}||\Pi_{S_{(k),(t)}} \epsilon||_n^2}{D_{k,t} ||\mu+\epsilon-\Pi_{V_{(k),(t)}}\mu-\Pi_{V_{(k),(t)}}\epsilon||_n^2}>x\right)\cap A_k\right).
\end{eqnarray*}

Since $A_k=\{\{X_{(1)},\dots,X_{(k)}\}=\{X_j,j\in J\}\}$, \\

$\mathbb{P}\left(\left(\dfrac{N_{k,t}||\Pi_{S_{(k),(t)}} \mu||_n^2+N_{k,t}||\Pi_{S_{(k),(t)}} \epsilon||_n^2}{D_{k,t} ||\mu+\epsilon-\Pi_{V_{(k),(t)}}\mu-\Pi_{V_{(k),(t)}}\epsilon||_n^2}>x\right)\cap A_k\right)$
\begin{eqnarray*}
\hspace{1cm}&=&\mathbb{P}\left(\left(\dfrac{N_{k,t}||\Pi_{S_{(k),(t)}} \epsilon||_n^2}{D_{k,t} ||\epsilon-\Pi_{V_{(k),(t)}}\epsilon||_n^2}>x\right)\cap A_k\right)\\
&\le&\mathbb{P}\left(\dfrac{N_{k,t}||\Pi_{S_{(k),(t)}} \epsilon||_n^2}{D_{k,t} ||\epsilon-\Pi_{V_{(k),(t)}}\epsilon||_n^2}>x\right).
\end{eqnarray*}
And by construction of $\Upsilon_{k,t}$,
$$\mathbb{P}\left(\dfrac{N_{k,t}||\Pi_{S_{(k),(t)}} \epsilon||_n^2}{D_{k,t} ||\epsilon-\Pi_{V_{(k),(t)}}\epsilon||_n^2}>x\right)\le\mathbb{P}\left(\Upsilon_{k,t}>x\right).$$
Thus $$\mathbb{P}\left((\tilde{U}_{D_{k,t},N_{k,t}}>x)\cap A_k\right)\le \mathbb{P}\left(\Upsilon_{k,t}>x\right).$$
\end{proof}

\begin{proof}[\textbf{Proof of Theorem \ref{th_3}}]
Let $k<k_0$ and $0<\gamma<1$. Denote $I=\{0,\dots,\lfloor log_2(k_0-k)\rfloor\}$. \\
From  the proof of Theorem \ref{th_2} (more precisely the condition $(cond_k)$), we have that 
 if the following condition is verified: \\
$ \exists t\in I$ such that
\begin{equation}
\label{cond1}
 \dfrac{1}{2} \text{inf}\{|| \Pi_{S}\mu||_n^2,S\in B_{2^t}\} \ge  \bar{\Upsilon}_{k,t}^{-1}(\alpha_{k,t})Q_{1-\gamma/2k_0}\dfrac{D_{k,t}}{N_{k,t}} + \dfrac{\sigma^2}{n}\bar{\chi}_{2^t}^{-1}\left(\dfrac{\gamma/2k_0}{|B_{2^t}|}\right),
\end{equation}
where $Q_{1-u}$ denote the $(1-u)$-quantile of the statistics $||Y-\Pi_{V_{(k),(t)}}Y||_n^2$ under the event $A_{k_0}$,\\
 then we have:
$$\mathbb{P}\left(\forall t\in I,||\Pi_{S_{(k),(t)}} Y||_n^2\le \bar{\Upsilon}_{k,t}^{-1}(\alpha_{k,t}) Q_{1-\gamma/2k_0}\dfrac{D_{k,t}}{N_{k,t}} \ \cap A_{k_0}\right)\le \gamma/2k_0.$$ 
Since \\
$\mathbb{P}\left(\forall t\in I, \tilde{U}_{D_{k,t},N_{k,t}}<\bar{\Upsilon}_{k,t}^{-1}(\alpha_{k,t})\ \cap A_{k_0}\right)$
\begin{eqnarray*}
\hspace{1cm}\le \underset{t\in I}{\text{inf }}\{ \mathbb{P}\left(\tilde{U}_{D_{k,t},N_{k,t}}<\bar{\Upsilon}_{k,t}^{-1}(\alpha_{k,t})\ \cap A_{k_0}\right)\}
\end{eqnarray*}
and since\\
$\mathbb{P}\left(\tilde{U}_{D_{k,t},N_{k,t}}<\bar{\Upsilon}_{k,t}^{-1}(\alpha_{k,t})\ \cap A_{k_0}\right)$
\begin{eqnarray*}
\hspace{1cm}&\le& \underset{\le \gamma/2k_0}{\underbrace{\mathbb{P}\left(||Y-\Pi_{V_{(k),(t)}}Y||_n^2>Q_{1-\gamma/2k_0} \ \cap A_{k_0}\right)}}\\
&&+\mathbb{P}\left(\dfrac{||\Pi_{S_{(k),(t)}} Y||_n^2}{D_{k,t}}\le \bar{\Upsilon}_{k,t}^{-1}(\alpha_{k,t}) \dfrac{Q_{1-\gamma/2k_0}}{N_{k,t}} \ \cap A_{k_0}\right)\\
&\le& \gamma/k_0,
\end{eqnarray*}
we have that the condition \eqref{cond1} implies that 
\begin{equation}
\label{impliesB}
\mathbb{P}\left(\forall t\in I, \tilde{U}_{D_{k,t},N_{k,t}}<\bar{\Upsilon}_{k,t}^{-1}(\alpha_{k,t})\ \cap A_{k_0}\right)\le \gamma/k_0.
\end{equation}
In the following, we give an upper bound of the right part in \eqref{cond1}. For this doing, we have to give an upper bound of $\bar{\Upsilon}_{k,t}^{-1}(\alpha_{k,t})$ and $Q_{1-\gamma/2k_0}$.\\

 Assume we are on the event $A_k$, then 
 \begin{eqnarray*}
 \Upsilon_{k,t}&=&\dfrac{N_{k,t} ||\Pi_{S_{(k),\sigma_1(t)}} Y||_n^2}{D_{k,t} ||Y-\Pi_{V_{(k),\sigma_1(t)}}Y||_n^2}\\
 &=&\dfrac{N_{k,t} ||\Pi_{S_{(k),\sigma_1(t)}} \epsilon||_n^2}{D_{k,t} ||Y-\Pi_{V_{(k)}}Y-\Pi_{S_{(k),\sigma_1(t)}}\epsilon||_n^2}.
 \end{eqnarray*}
As we are on the event $A_k$, the space $V_{(k)}$ is not a random space. Thus for any subspaces $S$ of dimension $D_{k,t}=2^t$, we have that $||\Pi_S Y||_n^2=||\Pi_S\epsilon||_n^2\sim \sigma^2\chi^2_{2^t}/n$ and we have that $ ||Y-\Pi_{V_{(k)}}Y-\Pi_{S}Y||_n^2=||\Pi_{(S\oplus V_{(k)})^{\bot}} \epsilon||_n^2\sim \sigma^2\chi^2_{n-(2^t+k)}/n$.\\ Hence  $\dfrac{N_{k,t} ||\Pi_{S} Y||_n^2}{D_{k,t} ||Y-\Pi_{V_{(k)}}Y-\Pi_{S}Y||_n^2}\sim F_{D_{k,t},N_{k,t}}$. Thus on the\\ event $A_k$,  $\Upsilon_{k,t}\le \text{sup}\{\dfrac{N_{k,t} ||\Pi_{S} \epsilon||_n^2}{D_{k,t} ||\epsilon-\Pi_{V_{(k)}+S}\epsilon||_n^2}, S \in G_{k,2^t}\}$,\\ where $G_{k,2^t}$ is defined by \eqref{Gkd}.\\
We deduce that the $(1-u)$-quantile of $\Upsilon_{k,t}$ is lower that\\ $\bar{F}_{D_{k,t},N_{k,t}}^{-1}\left(u/|G_{k,2^t}|\right)$. Indeed: \\
$\mathbb{P}\left(\Upsilon_{k,t}>\bar{F}_{D_{k,t},N_{k,t}}^{-1}\left(u/|G_{k,2^t}|\right)\right)$
\begin{eqnarray*}
 &\le &\mathbb{P}\left[  \text{sup}\{\dfrac{N_{k,t} ||\Pi_{S} \epsilon||_n^2}{D_{k,t} ||\epsilon-\Pi_{V_{(k)}+S}\epsilon||_n^2}, S \in G_{k,2^t}\}\right.\\
&&\hspace{1cm} > \left.\bar{F}_{D_{k,t},N_{k,t}}^{-1}\left(u/|G_{k,2^t}|\right)\right]\\
&\le&\sum_{S \in G_{k,2^t}} \mathbb{P}\left(\dfrac{N_{k,t} ||\Pi_{S} \epsilon||_n^2}{D_{k,t} ||\epsilon-\Pi_{V_{(k)}+S}\epsilon||_n^2}>\bar{F}_{D_{k,t},N_{k,t}}^{-1}\left(u/|G_{k,2^t}|\right)\right)\\
&\le&|G_{k,2^t}| \dfrac{u}{|G_{k,2^t}|}\le u.
\end{eqnarray*}
\citet{Baraud:2003} gave an upper bound of $\bar{F}_{D,N}^{-1}\left(u\right)$, for $0<D$, $0<N$ and $0<u$:
 \begin{eqnarray*}
 D\bar{F}_{D,N}^{-1}(u) &\le& D+2 \sqrt{D\left(1+\dfrac{D}{N}\right)log\left(\dfrac{1}{u}\right)}\\
 &&+\left(1+2\dfrac{D}{N}\right)\dfrac{N}{2}\left[exp\left(\dfrac{4}{N}log\left(\dfrac{1}{u}\right)\right)-1\right].
 \end{eqnarray*}
Since $\text{exp}(u)-1\le u \text{ exp}(u)$ for any $u>0$,  $\sqrt{u+v}\le \sqrt{u}+\sqrt{v}$ for all $u>0,v>0$ and since  $\alpha_{k,t}\ge \alpha/|\mathcal{T}_k|$, we derive that:
\begin{eqnarray*}
2^t\bar{\Upsilon}_{k,t}^{-1}(\alpha_{k,t})  &\le& 2^t\left[ 1+\Lambda_3(k,t)log\left(\dfrac{e(p-k)}{2^t}\right)\right]\\
&+& 2\left[\sqrt{2^t\left(1+\dfrac{2^t}{N_{k,t}}\right) log\left(\dfrac{|\mathcal{T}_k|}{\alpha}\right)}+\dfrac{\Lambda_2(k,t)}{2}log\left(\dfrac{|\mathcal{T}_k|}{\alpha}\right)\right],
  \end{eqnarray*}
where $\Lambda_1(k,t)=\sqrt{1+\dfrac{D_{k,t}}{N_{k,t}}}$, $\Lambda_2(k,t)=\left(1+2\dfrac{D_{k,t}}{N_{k,t}}\right)M$ and $\Lambda_3(k,t)=2\Lambda_1(k,t)+\Lambda_2(k,t)$ with $L_t=\text{log}(|\mathcal{T}_k|/\alpha)$,
$m_t=\text{exp}(4L_t/N_{k,t})$, $ m_p=\text{exp}\left(\dfrac{4D_{k,t}}{N_{k,t}}log\left(\dfrac{e(p-k)}{2^t}\right)\right)$,
$M=2m_tm_p$.\\
Since $\sqrt{ab}+mb\le a/2+(m+1/2)b$ holds for any positive numbers $a,b,m$, we obtain that:

\begin{eqnarray}
  \label{gam1}
 2^t\bar{\Upsilon}_{k,t}^{-1}(\alpha_{k,t})   &\le& 2^t\left[ 1+\Lambda_1^2(k,t)+\Lambda_3(k,t)log\left(\dfrac{e(p-k)}{2^t}\right)\right]\\
 &+&\left(1+\Lambda_2(k,t)\right)log\left(\dfrac{|\mathcal{T}_k|}{\alpha}\right).
  \label{gam2}
\end{eqnarray}

We have now to find an upper bound of $Q_{1-\gamma/2k_0}$.\\
$Q_{1-\gamma/2k_0}$ is defined by $\mathbb{P}\left(||Y-\Pi_{V_{(k),(t)}}Y||_n^2>Q_{1-\gamma/2k_0} \ \cap A_{k_0}\right)\le \gamma/2k_0$.\\
We always have that:
$||Y-\Pi_{V_{(k),(t)}}Y||_n^2 \le ||\mu||^2_n+||\epsilon||_n^2$.
Thus $\forall \ 0<u<1$, the $(1-u)$-quantile of $||Y-\Pi_{V_{(k),(t)}}Y||_n^2 $ is lower than the $(1-u)$-quantile of $||\mu||^2_n+||\epsilon||_n^2$.\\
As $||\epsilon||_n^2\sim \sigma^2\chi_n^2/n$, we can use the equation \eqref{chi2} for $x_u=log(2k_0/\gamma)$ and we obtain that $\bar{\chi}_n^{-1}(\gamma/2k_0)\le n+2\sqrt{nx_u}+2x_u$.\\
Therefore 
\begin{equation}
\label{Qgam1}
Q_{1-\gamma/2k_0}\le ||\mu||^2_n+\sigma^2\dfrac{n+2\sqrt{nx_u}+2x_u}{n}
\end{equation}
and as $1+2\sqrt{u}+2u\le 2+3u$, we get 
\begin{equation}
\label{Qgam}
Q_{1-\gamma/2k_0}\le ||\mu||^2_n+\sigma^2\left(2+\dfrac{3}{n}log\left(\dfrac{2k_0}{\gamma}\right)\right).
\end{equation}

Combining  \eqref{gam1}, \eqref{Qgam} in \eqref{cond1} and using that 
\begin{eqnarray*}
\bar{\chi}_{2^t}^{-1}\left(\dfrac{\gamma/2k_0}{|B_{2^t}|}\right)&\le& 2^t\left[5 +4log\left(\dfrac{k_0}{2^t}\right)\right]+2\left[\sqrt{2^t log(2k_0/\gamma)}+log(2k_0/\gamma)\right]\\
&\le & 2^t\left[6 +4log\left(\dfrac{k_0}{2^t}\right)\right]+3log(2k_0/\gamma) ,
\end{eqnarray*}
 we obtain the following condition:\\
$(R_{3,k}):\ \exists t\in I$ such that\\
$ \dfrac{1}{2} \text{inf}\{|| \Pi_{S}\mu||_n^2,S\in B_{2^t}\}$
\begin{eqnarray*}
 &\ge&   \dfrac{D_{k,t}\bar{F}_{D_{k,t},N_{k,t}}^{-1}(\alpha_{k,t}/|G_{2^t}|)}{N_{k,t}} \left[||\mu||^2_n+\sigma^2\left(2+\dfrac{3}{n}log\left(\dfrac{2k_0}{\gamma}\right)\right)\right]\\
 & & +\dfrac{\sigma^2}{n}\left[2^t\left(6 +4log\left(\dfrac{k_0}{2^t}\right)\right)+3log\left(\dfrac{2k_0}{\gamma}\right)\right]\\
&\ge& \dfrac{A(k,t)}{N_{k,t}}
\left[||\mu||^2_n+\sigma^2\left(2+\dfrac{3}{n}log\left(\dfrac{2k_0}{\gamma}\right)\right)\right]\\
 && +\dfrac{\sigma^2}{n}\left[2^t\left(6 +4log\left(\dfrac{k_0}{2^t}\right)\right)+3log\left(\dfrac{2k_0}{\gamma}\right)\right],
 \end{eqnarray*}
 where 
 $A(k,t)= 2^t\left[ 2+\dfrac{2^{t}}{N_{k,t}}+\Lambda_3(k,t)log\left(\dfrac{e(p-k)}{2^t}\right)\right]
 +\left(1+\Lambda_2(k,t)\right)log\left(\dfrac{|\mathcal{T}_k|}{\alpha}\right)
$.\\\\
 The condition $(R_{3,k})$ leads to \eqref{impliesB} and thus 
\begin{eqnarray*}
\mathbb{P}_{\mu}(\hat{J}\neq J )& \le&\mathbb{P}_{\mu}(\hat{J}\neq J \ \cap  A_{k_0})+\mathbb{P}( A_{k_0}^c) \\
&\le &\left(\sum_{j=0}^{k_0-1} \mathbb{P}_{\mu}(\mathring{k}_{B}=j \ \cap  A_{k_0})+\mathbb{P}_{\mu}(\mathring{k}_{B}>k_0 \ \cap A_{k_0})\right)+\mathbb{P}( A_{k_0}^c) \\
&\le& k_0  \gamma/k_0 + \alpha + \delta.
\end{eqnarray*}
And then \eqref{non3} is proved.
\end{proof}

 \begin{proof}[\textbf{Proof of Remark \ref{rem3}}]
 In the following, $C(a,b)$ denote a constant depending on the parameters $a$ and $b$.
 Under the assumption that  $2^t\le (n-k)/2$ and since $\forall \ x\ge 2, \dfrac{log(x)}{x}\le 1$ we have that:
 $$\dfrac{D_{k,t}}{N_{k,t}}log\left(\dfrac{p-k}{D_{k,t}}\right)\le \dfrac{2^t}{n-k-2^t}log\left(\dfrac{n-k}{2^t}\right)\le 2\dfrac{2^t}{n-k}log\left(\dfrac{n-k}{2^t}\right)\le 2.$$

Moreover the ratio $D_{k,t}/N_{k,t}$ is bounded by $1$,  thus $log(m_p)\le 4\dfrac{D_{k,t}}{N_{k,t}}+4\dfrac{D_{k,t}}{N_{k,t}}log\left(\dfrac{p-k}{D_{k,t}}\right)\le 12 $.\\
 As the ratio $4L_{k,t}/N_{k,t}$ is bounded by $C'(\alpha)$ and since \\$M\le 2\text{exp}(C'(\alpha)) \text{exp}(12)$, we have that $M$ is bounded by $C''(\alpha)$. Thus 
 $\Lambda_1(k,t)\le \sqrt{2}$, $\Lambda_2(k,t)\le 3C''(\alpha)$ and $\Lambda_3(k,t)\le 2\sqrt{2}+3C''(\alpha)$.\\
 We obtain under the condition $log(p-k)>1$ that $A(k,t)\le 2^t C(\alpha) log(p-k)$.\\
We also have that
$\left[||\mu||^2_n+\sigma^2\left(2+\dfrac{3}{n}log\left(\dfrac{2k_0}{\gamma}\right)\right)\right]\le C(||\mu||_n,\gamma,\sigma)$ \\since $log(k_0)/n\le 1$,
and that\\ $2^t\left(6 +4log\left(\dfrac{k_0}{2^t}\right)\right)+3log\left(\dfrac{2k_0}{\gamma}\right)\le 2^t\left[6 +4log\left(\dfrac{k_0}{2^t}\right)+3log\left(\dfrac{2k_0}{\gamma}\right)\right]\le 2^tC(\gamma)log(k_0)$.\\\\

We finally obtain equation \eqref{simp}.
 \end{proof}

\begin{proof}[\textbf{Proof of Remark \ref{ortho}}]
The differences between the two conditions $(R_{3,k})$ and $(R_{3bis,k})$ lie in the fact that \\ $\text{inf}\{|| \Pi_{S}\mu||_n^2,S\in B_{2^t}\}=\sum_{j=1}^{2^t} \beta^2_{\sigma_2(j)}$ and that the upper bound of $Q_{1-\gamma/2k_0}$ is modified, where $Q_{1-\gamma/2k_0}$ is defined by \\$\mathbb{P}\left(||Y-\Pi_{V_{(k),(t)}}Y||_n^2>Q_{1-\gamma/2k_0} \ \cap A_{k_0}\right)\le \gamma/2k_0$. \\Indeed, on the event ${A}_{k_0}$ we have that
$||Y-\Pi_{V_{(k),(t)}}Y||_n^2 \le \sum_{j=k+2^t}^{j=k_0}\beta^2_{\sigma_2(j)} +||\epsilon||_n^2$, where $\sigma_2$ is defined such that $|\beta_{\sigma_2(1)}|\le ... \le |\beta_{\sigma_2(k_0)}|$. We get from there the condition $(R_{3bis,k})$.
\end{proof}

\end{document}